\documentclass[final,3p]{elsarticle}
\usepackage{color}
\usepackage{amssymb}
\usepackage{amsmath}
\usepackage{booktabs}
\usepackage{multirow}
\usepackage{tabularx}
\usepackage{tabu}
\usepackage{stmaryrd}
\usepackage{epstopdf}
\usepackage{subfigure} 
\usepackage{graphicx}
\usepackage{caption}
\usepackage{amsthm} 
\usepackage{amssymb} 
\usepackage{mathrsfs}
\newcommand{\ds}{\displaystyle}

\newcommand{\mbf}{\mathbf}

\numberwithin{equation}{section}
\newtheorem{thm}{Theorem}[section]

\newtheorem{rem}[thm]{Remark}
\allowdisplaybreaks[4]

\biboptions{}
\journal{*}
\begin{document}
	
\begin{frontmatter}
\title{Energy-preserving fully-discrete schemes for nonlinear stochastic wave equations with multiplicative noise}

%---------作者----------------------------------------------------------------------------------
\author[ad1,ad2]{Jialin Hong}
\ead{	hjl@lsec.cc.ac.cn} 

\author[ad1]{Baohui Hou\corref{cof1}}
\ead{houbaohui@lsec.cc.ac.cn} 

\author[ad1]{Liying Sun}
\ead{liyingsun@lsec.cc.ac.cn}

\cortext[cof1]{Corresponding author}
\address[ad1]{Institute of Computational Mathematics and Scientific/Engineering Computing, Academy of Mathematics and Systems Science, Chinese Academy of Sciences, Beijing 100190, China}
\address[ad2]{School of Mathematical Sciences, University of Chinese Academy of Sciences, Beijing 100049, China}
%---------------abstract----------------------------------------------------------------------------		
\begin{abstract}	
In this paper, we focus on constructing numerical schemes preserving the averaged energy evolution law for nonlinear stochastic wave equations  driven by multiplicative noise. We first apply the compact finite difference method and the interior penalty discontinuous Galerkin finite element method to discretize  space variable and present two semi-discrete schemes, respectively.  Then we make use of  the discrete gradient method and the Pad\'e approximation to propose efficient fully-discrete schemes. These semi-discrete and fully-discrete schemes are proved to preserve the discrete averaged energy evolution law.  
In particular, we also prove that the proposed fully-discrete schemes exactly inherit the averaged  energy evolution law almost surely if the considered model is driven by additive noise. Numerical experiments are given to confirm theoretical findings.

\end{abstract}		
		
		
\end{frontmatter}
%-----------section 1----Introduction---------------------------------
\section{Introduction}
The nonlinear stochastic wave equation plays an important role in a wide range of applications in the field of  engineering,  science, etc.,  and is commonly used to describe a variety of physical processes, such as the motion of a strand of DNA in a liquid, the motion of shock waves on the surface of the sun, the dynamics of the primary current density vector field within the grey matter of the human brain and the sound propagation in the sea and so on (see e.g., \cite{Banjai,Chow, Cohen-Larsson, Dalang} and references therein). 
In this paper, we consider the following nonlinear stochastic wave equation driven by multiplicative noise
\begin{equation} \label{eq1.1}\left\{\begin{aligned}
&du = vdt,  &&(x, t) \in (a, b) \times (0,T],\\
&dv =  \Delta u dt - f(u) dt +g(u)dW(t), &&(x, t) \in (a, b) \times (0,T],\\
& u(x,0) = u_0(x), v(x,0) = v_0(x), &&x \in [a, b],
\end{aligned} \right. \end{equation}
 where $u_0$ and $v_0$ are real-valued deterministic functions, and $\Delta$ is the Laplace operator with Dirichlet or periodic boundary condition.  Here, $\{W(t)\}_{t\geq 0}$ is an $\mathscr H$-valued $Q$-Wiener process with respect to a normal filtration $\{\mathcal{F}_t\}_{t\geq 0}$ on a filtered probability space $(\Omega, \mathcal{F},\{\mathcal{F}_t\}_{t\geq 0}, \mathbb{P})$ and has the form of
 $$W(t)=\sum_{k=1}^{\infty} \sqrt{q_k}e_k\beta_k(t),$$ 
 where $\mathscr H  = L^2([a,b];\mathbb{R}) $,  $\{(q_k, e_k)\}_{k=1}^{\infty}$ are eigenpairs of  symmetric, positive definite and finite trace operator $Q$ with orthonormal eigenvectors and $\{\beta_k(t)\}_{k=1}^{\infty}$ is a sequence of real-valued mutually independent standard Brownian motions.
Assume that $f: \mathscr H \rightarrow  \mathscr H$ and $ g(u): \mathscr H \rightarrow \mathscr{L}_2^0$ satisfy the global Lipschitz condition, where $\mathscr{L}_2^0 = \mathscr{L}_2(\mathscr H, \mathscr H_0)$ denotes the  separable Hilbert space of Hilbert-Schmidt operators from $\mathscr H$ to $\mathscr H_0$ with $\mathscr H_0 = Q^{\frac{1}{2}}(\mathscr H)$. For the well-posedness of \eqref{eq1.1}, we refer readers to \cite {Anton, Wang}.  It has been shown that the averaged energy of the nonlinear  stochastic wave equation \eqref{eq1.1} satisfies the following identity
\begin{equation}\label{eqE}
\mathbb E[H(u(t),v(t))] = \mathbb E[H(u(0),v(0))] +\frac{1}{2}\mathbb E\left[\int_{0}^{t} Tr\left((g(u) Q^\frac{1}{2})(g(u) Q^\frac{1}{2})^{*}\right)ds\right],
\end{equation}
where
\begin{equation}\label{eq1.2}
H(u,v)=\frac{1}{2}\int_{a}^{b}\left(v^2+\left(\frac{\partial u}{\partial x}\right)^2\right)dx+\int_{a}^{b} \widetilde f\left(u\right)dx,
\end{equation}
 and $f=\frac{\partial  \widetilde f}{\partial u}$.
In particular, when $g(u)=\lambda$ with $\lambda$ being a positive constant, \eqref{eq1.1} becomes a nonlinear stochastic wave equation driven by additive noise. It is shown in \eqref{eqE} that the averaged energy increases linearly with respect to the evolution of time with a growth rate $\frac{\lambda^2}{2}Tr(Q)$. \eqref{eqE} is also related to the energy equation,  which is a tool that can be used to analyze the existence or nonexistence of solutions to nonlinear stochastic wave equations (see \cite{Anton}).

It is well-known that numerical simulations are often used to understand the behavior of the solutions for stochastic partial differential equations. Especially, numerical schemes preserving energy evolution law often yield physically correct results and numerical stability in computation. Much effort has been devoted to the numerical method of stochastic wave equations (see e.g.,  \cite{Anton, Banjai, Cao,Cohen-Larsson,Cohen,Cui, Kovacs, Quer, Walsh} and references therein).
For instance, \cite{Cohen-Larsson} proves that the  fully-discrete scheme based on the finite element method and the stochastic trigonometric scheme preserves the linear growth of the averaged energy for the linear stochastic wave equation with additive noise.
With regard to the nonlinear stochastic wave equation with additive noise, \cite{Banjai} shows that the discontinuous Galerkin finite element method satisfies the trace formula. \cite{Cui} proposes a full discretization by adopting the spectral Galerkin method and the averaged vector field method preserving the averaged energy evolution law of the stochastic cubic wave equation with additive noise. In \cite{Anton}, the authors utilize the standard linear finite element approximation and a stochastic trigonometric method to propose a fully-discrete scheme for the nonlinear stochastic wave equation driven by multiplicative noise, and  prove that this numerical scheme satisfies an almost trace formula for the case of the additive noise. To the best of knowledge, there has not been a fully-discrete scheme which can exactly preserve the averaged energy evolution law of nonlinear stochastic wave equations driven by multiplicative noise. 

In this paper, we aim to propose fully-discrete schemes preserving averaged energy evolution law for nonlinear stochastic wave equations driven by multiplicative noise \eqref{eq1.1}. 
Motivated by the fact that the compact finite difference method enjoys the flexibility to handle the boundary conditions and could achieve high-order accuracy with smaller stencils and the fact that the interior penalty discontinuous Galerkin finite element method is particularly suitable to deal with complex computational domains and is easy to design high-order approximations, we apply these two numerical schemes to discretizing \eqref{eq1.1} and obtain two semi-discrete schemes in space. 
We prove that the resulting finite dimensional stochastic differential equations preserve the discrete version of the averaged energy evolution law of the original system. 
Besides numerical schemes in space, the energy-preserving fully-discrete scheme also depends on the numerical schemes in time,
which confronts the difficulty brought by the treatment of the time approximation on both drift and diffusion coefficients. 
For example, if we make use of the discrete gradient method to the nonlinear drift term and the explicit Euler method to the diffusion term for semi-discrete scheme, which is given by means of the compact finite difference in spatial direction,  of the nonlinear stochastic wave equation \eqref{eq1.1}, the averaged energy of the obtained fully-discrete scheme is more $\frac{1}{2}h\Delta t \mathbb E\left[(\frac{1}{2} \mbf {A^{-1}D}(\mbf U^{n+1}+\mbf U^n)-\overline{\nabla}_U\mbf F(\mbf U^{n+1},\mbf U^n))G(\mbf U^n)\mbf E \boldsymbol{\Lambda} \Delta\boldsymbol {\beta}_n\right]$ than that of the original stochastic system (the concrete meanings of symbols are given below).  
To overcome this difficulty, we discretize the solution given by the constant variation method of aforementioned two semi-discrete schemes based on the discrete gradient method and the Pad\'e approximation. We prove that the proposed fully-discrete schemes admit the discrete averaged energy evolution law, which is consistent with the averaged energy evolution law of \eqref{eq1.1}. 
We would like to mention that the proposed fully-discrete schemes are flexible due to the flexibility of the Pad\'e approximation. 
Finally, numerical experiments confirm the theoretical analysis results.

%% 剩余文章结构
The paper is organized as follows.  In Section 2,  we employ the compact finite difference method and the interior penalty discontinuous Galerkin finite element method to obtain two semi-discrete schemes, and prove that both of them possess discrete averaged energy evolution laws.  Then the discrete gradient method and the Pad\'e approximation are ulitized to construct fully-discrete schemes,  and the conservation of  discrete averaged energy evolution laws is proved in Section 3.  Numerical experiments are carried out in Section 4.
Finally, the conclusion is given in Section 5.

%-----------section 2------------------------------------------------------------
\section{Energy-preserving semi-discrete schemes in space}
In this section, we investigate the compact finite difference method and the interior penalty discontinuous Galerkin finite element method for the nonlinear stochastic wave equation \eqref{eq1.1}, respectively. Then we prove that these semi-discrete schemes preserve the averaged energy evolution laws for \eqref{eq1.1} below. Denote the uniform partition of $[a, b]$ by $\{x_i\}_{0\leq i\leq M}$ with $x_i =  a+ i h$, where $M$ is a positive integer, and $h := (b-a)/M$ denotes the spatial step size.  

%------------------CFD半离散-----------------------------------
\subsection{Semi-discrete scheme via the compact finite difference method}
The compact finite difference method, as the high accuracy discrete method utilizing smaller stencils and treating boundary conditions easily, has been extensively studied in \cite{Britt, Deng, Moghaderi}. 
We first apply the compact finite difference method to discretizing the nonlinear stochastic wave equation \eqref{eq1.1} with Dirichlet boundary condition.
Now we introduce the discrete grid $\Omega_h = \{ x_i |~ 1 \leq i \leq M-1\}$, and difference operators $\delta_x^2$ and $\mathcal A$ defined as
$$\quad \delta_x^2u_{i} = \frac{u_{i-1}-2u_{i}+u_{i+1}}{h^2}, \qquad \mathcal A u_{i} = \left(1+\frac{h^2}{12}\delta_x^2\right)u_{i}$$
for any grid function $u = \{ u_i |~ x_i \in \Omega_h, u_0=u_M=0\}$.
From \cite{Hirsh,Moghaderi}, it can be verified that the compact finite difference operator $\mathcal A^{-1} \delta_x^2$ satisfies 
\begin{equation*}
\mathcal A^{-1} \delta_x^2u_{i} = \frac{\partial^2 u}{\partial x^2}(x_i) +O(h^4), \quad 1\leq i\leq M-1.
\end{equation*}
Based on the above operator, the compact finite difference method for \eqref{eq1.1} with Dirichlet boundary  condition reads
 \begin{equation} \label{eq2.1}\left\{\begin{aligned}
dU_i &=V_idt,\\
dV_i &= \mathcal A^{-1} \delta_x^2U_{i}dt -f(U_{i})dt 
+ g(U_{i})\sum_{k=1}^{P}\sqrt{q_k}e_k(x_i)d\beta_k(t),
\end{aligned} \right. \end{equation} 
where $U_i, V_i$ are the approximations of $u_i= u(x_i,t), v_i = v(x_i,t)$, respectively, for $1\leq i \leq M-1 $. 
Denote 
\begin{equation*}\begin{aligned}
&\mbf U = \left(U_1, U_2, \dots,U_{M-1}\right)^\top,\quad&&\mbf V = \left(V_1, V_2, \dots,V_{M-1}\right)^\top,\\
&\boldsymbol{\Lambda}= {\rm{diag} }\{ \sqrt{q_1}, \sqrt{q_2}, \dots,\sqrt{q_P} \},\quad &&\boldsymbol {\beta}=\left(\beta_1, \beta_2, \dots,\beta_P\right)^\top,\\
&F(\mbf U)=\left(f(U_{1}), f(U_{2}), \dots,f(U_{M-1})\right)^\top,\quad &&G(\mbf U)={\rm{diag} }\{ g(U_{1}), g(U_{2}), \dots,g(U_{M-1})\}
\end{aligned} \end{equation*}	
 and 
\begin{equation*}
\mbf E={
	\left[ \begin{array}{cccc}
	e_1(x_1)&e_2(x_1)&\cdots&e_P(x_1)\\
	e_1(x_2)&e_2(x_2)&\cdots&e_P(x_2)\\
	\vdots&\vdots&&\vdots\\
	e_1(x_{M-1})&e_2(x_{M-1})&\cdots&e_P(x_{M-1})\\
	\end{array}
	\right ]_{(M-1) \times P}}.
\end{equation*}	
Define matrices $\mbf{D}$ and $\mbf{A}$ of order $(M-1)\times (M-1)$ as
\begin{equation}\label{eqDA}
\mbf{D}=\frac{1}{h^2}{
	\left[ \begin{array}{cccccccc}
	-2&1&0&0&\cdots&0&0&0\\[0.05in]
	1&-2&1&0&\cdots&0&0&0\\
	\vdots&\vdots&\vdots&\vdots&&\vdots&\vdots&\vdots\\
	0&0&0&0&\cdots&1&-2&1\\[0.05in]
	0&0&0&0&\cdots&0&1&-2
	\end{array}
	\right ]},\quad
\mbf{A}={
	\left[ \begin{array}{cccccccc}
	\frac{10}{12}&\frac{1}{12}&0&0&\cdots&0&0&0\\[0.05in]
	\frac{1}{12}&\frac{10}{12}&\frac{1}{12}&0&\cdots&0&0&0\\
	\vdots&\vdots&\vdots&\vdots&&\vdots&\vdots&\vdots\\
	0&0&0&0&\cdots&\frac{1}{12}&\frac{10}{12}&\frac{1}{12}\\[0.05in]
	0&0&0&0&\cdots&0&\frac{1}{12}&\frac{10}{12}
	\end{array}
	\right ]},
\end{equation}	
which are arising from the operator $\delta_x^2$ and $ \mathcal A$, respectively.  In addition, let $\mbf{A^{-1}D}$ be the matrix associated to operator $\mathcal A^{-1} \delta_x^2$.
Then we can rewrite \eqref{eq2.1} into 
\begin{equation} \label{eq2.2}\left\{\begin{aligned}
d\mbf U &=\mbf Vdt,\\
d\mbf V&= \mbf{A^{-1}D}\mbf Udt -F(\mbf U)dt 
+ G(\mbf U)\mbf E \boldsymbol{\Lambda}d \boldsymbol {\beta}.\\
\end{aligned} \right. \end{equation} 

The semi-discrete scheme \eqref{eq2.2} preserves the discrete version of stochastic energy evolution law, which is stated in the following theorem.
\begin{thm} \label{thm2}
The averaged energy for semi-discrete scheme \eqref{eq2.2} has the following evolutionary relationship
\begin{equation}
\mathbb E\left[\bar{H}(\mbf U(t),\mbf V(t))\right]=\mathbb E\left[\bar{H}(\mbf U(0),\mbf V(0))\right]
+\frac{h}{2}\mathbb E\left[\int_0^tTr(G^2(\mbf U)\mbf E\boldsymbol{\Lambda}^2\mbf E^\top)ds\right],
\end{equation}
where 
\begin{equation*}
\bar{H}(\mbf U(t),\mbf V(t))=\frac{1}{2}\sum_{i=1}^{M-1}hV_i^2(t)-\frac{1}{2}\sum_{i=1}^{M-1}hU_i(t) \mathcal A^{-1} \delta_x^2U_i(t)+\sum_{i=1}^{M-1}h\widetilde f(U_i(t))
\end{equation*}
is the discrete energy with respect to \eqref{eq1.2}.
\end{thm}

\noindent{\textbf{Proof. }}
Based on It\^o's formula for $\frac{1}{2}\sum_{i=1}^{M-1}hV_i^2$, we obtain
\begin{equation}\begin{aligned}\label{eq2.3}
\frac{1}{2}\sum_{i=1}^{M-1}hV_i^2(t)&=\frac{1}{2}\sum_{i=1}^{M-1}hV_i^2(0) 
+ \int_{0}^{t}\sum_{i=1}^{M-1}hV_i(s)g(U_i(s))\sum_{k=1}^{P}\sqrt{q_k}e_k(x_i)d\beta_k(s)\\
&\quad+\int_{0}^{t}\sum_{i=1}^{M-1}hV_i(s) (\mathcal A^{-1} \delta_x^2U_i(s))ds
-\int_{0}^{t}\sum_{i=1}^{M-1}hV_i(s) f(U_i(s))ds\\
&\quad+\frac{1}{2}\int_{0}^{t} h\sum_{i=1}^{M-1}\sum_{k=1}^{P}g^2(U_i(s))q_ke_k^2(x_i)ds.
\end{aligned}\end{equation}
Direct calculations yield
\begin{equation}\begin{aligned}\label{eq2.4}
\int_{0}^{t}\sum_{i=1}^{M-1}hV_i(s) (\mathcal A^{-1} \delta_x^2U_i(s))ds
&=\int_{0}^{t}\sum_{i=1}^{M-1}\frac{h}{2}d\left(U_i(s) \mathcal A^{-1} \delta_x^2U_i(s)\right)\\
&=\sum_{i=1}^{M-1}\frac{h}{2}\left(U_i(t) \mathcal A^{-1} \delta_x^2U_i(t)-U_i(0) \mathcal A^{-1} \delta_x^2U_i(0)\right),
\end{aligned}\end{equation}
\begin{equation}\label{eq2.5}
-\int_{0}^{t}\sum_{i=1}^{M-1}hV_i(s) f(U_i(s))ds
=-\int_{0}^{t}\sum_{i=1}^{M-1}hd\widetilde f(U_i(s)) =-\sum_{i=1}^{M-1}h\left(\widetilde f(U_i(t)) -\widetilde f(U_i(0) )\right).
\end{equation}
Combining  the definition of $\bar{H}(\mbf U(t),\mbf V(t))$ as follows
\begin{equation}
\bar{H}(\mbf U(t),\mbf V(t))=\frac{1}{2}\sum_{i=1}^{M-1}hV_i^2(t)-\frac{1}{2}\sum_{i=1}^{M-1}hU_i(t) \mathcal A^{-1} \delta_x^2U_i(t)+\sum_{i=1}^{M-1}h\widetilde f(U_i(t))
\end{equation}
and \eqref{eq2.4}, \eqref{eq2.5},  then \eqref{eq2.3} can be rewritten as 
\begin{equation}\begin{aligned}
\bar{H}(\mbf U(t),\mbf V(t))&=\bar{H}(\mbf U(0),\mbf V(0))+ \int_{0}^{t}\sum_{i=1}^{M-1}hV_i(s)g(U_i(s))\sum_{k=1}^{P}\sqrt{q_k}e_k(x_i)d\beta_k(s)\\
&\quad+\frac{h}{2}\int_{0}^{t}\sum_{i=1}^{M-1}\sum_{k=1}^{P}g^2(U_i(s))q_ke_k^2(x_i)ds.
\end{aligned}\end{equation}
Taking expectation, we derive
\begin{equation*}
\mathbb E\left[\bar{H}(\mbf U(t),\mbf V(t))\right]=\mathbb E\left[\bar{H}(\mbf U(0),\mbf V(0))\right]
+\frac{h}{2}\mathbb E\left[\int_0^tTr(G^2(\mbf U)\mbf E\boldsymbol{\Lambda}^2\mbf E^\top)ds\right],
\end{equation*}
which completes the proof. 
\hfill $\square$\\

\begin{rem}
 If we consider the nonlinear  stochastic wave equation \eqref{eq1.1} with periodic boundary condition,  the corresponding matrices associated to operator $\delta_x^2$ and $ \mathcal A$ are 
 \begin{equation*}
\mbf{D}_{p}=\frac{1}{h^2}{
	\left[ \begin{array}{cccccccc}
	-2&1&0&0&\cdots&0&0&1\\[0.05in]
	1&-2&1&0&\cdots&0&0&0\\
	\vdots&\vdots&\vdots&\vdots&&\vdots&\vdots&\vdots\\
	0&0&0&0&\cdots&1&-2&1\\[0.05in]
	1&0&0&0&\cdots&0&1&-2
	\end{array}
	\right ]},\quad
\mbf{A}_{p}={
	\left[ \begin{array}{cccccccc}
	\frac{10}{12}&\frac{1}{12}&0&0&\cdots&0&0&\frac{1}{12}\\[0.05in]
	\frac{1}{12}&\frac{10}{12}&\frac{1}{12}&0&\cdots&0&0&0\\
	\vdots&\vdots&\vdots&\vdots&&\vdots&\vdots&\vdots\\
	0&0&0&0&\cdots&\frac{1}{12}&\frac{10}{12}&\frac{1}{12}\\[0.05in]
	\frac{1}{12}&0&0&0&\cdots&0&\frac{1}{12}&\frac{10}{12}
	\end{array}
	\right ]},
\end{equation*}	
where $\mbf{D}_{p}$ and $\mbf{A}_{p}$ are matrices of order $M\times M$.
Also, the corresponding semi-discrete scheme discretized by the compact finite difference method is 
\begin{equation} \label{Peq1}\left\{\begin{aligned}
d\mbf U &=\mbf Vdt,\\
d\mbf V&= \mbf{A}_p^{-1}\mbf{D}_p\mbf Udt -F(\mbf U)dt 
+ G(\mbf U)\mbf E_p \boldsymbol{\Lambda}d \boldsymbol {\beta},\\
\end{aligned} \right. \end{equation} 
where $\mbf U= \left(U_1, U_2, \dots,U_M\right)^\top$, $\mbf V=\left(V_1, V_2, \dots,V_M\right)^\top$, $F(\mbf U)=\left(f(U_{1}), f(U_{2}),\dots,f(U_{M}) \right)^\top$, $G(\mbf U) = {\rm{diag} }(g(U_{1}), g(U_{2}),\dots,g(U_{M}))$ and $\mbf E_p = (e_i(x_k))$ for $ i \in \{1,2,\dots, P\}, k\in\{1,2,\dots, M\}$.
Similarly,  the semi-discrete scheme \eqref{Peq1}  of the nonlinear stochastic wave equation with periodic boundary condition possesses the averaged energy evolution law, i.e.,
\begin{equation*}
\mathbb E\left[\bar{\bar{H}}(\mbf U(t),\mbf V(t))\right]=\mathbb E\left[\bar{\bar{H}}(\mbf U(0),\mbf V(0))\right]
+\frac{h}{2}\mathbb E\left[\int_0^tTr(G^2(\mbf U)\mbf E_p\boldsymbol{\Lambda}^2\mbf E_p^\top)ds\right],
\end{equation*}
where 
\begin{equation*}
\bar{\bar{H}}(\mbf U(t),\mbf V(t))=\frac{1}{2}\sum_{i=1}^{M}hV_i^2(t)-\frac{1}{2}h\mbf U(t)^\top\mbf{A}_p^{-1}\mbf{D}_p\mbf U+\sum_{i=1}^{M}h\widetilde f(U_i(t)).
\end{equation*}
\end{rem}

%--------DG半离散----------------------
\subsection{Semi-discrete scheme via the interior penalty discontinuous Galerkin finite element method}
The discontinuous Galerkin finite element method is flexible to deal with the complex computational domain and is easy to construct locally high-order approximations, which has been extensively studied in  \cite{Arnold, Banjai, Li}.
 In this subsection, we discretize  \eqref{eq1.1} with Dirichlet boundary condition by using the interior penalty discontinuous Galerkin finite element method in space.  Here, we use the same uniform mesh division of $[a, b]$ as in Subsection 2.1, denote $I_i=(x_i, x_{i+1})$, and define the discontinuous polynomial space
as follows 
\begin{equation*}
V_h=\{u\in L^2([a, b]): u|_{I_i} \in \mathcal{P}^k, ~u(a)=u(b)=0 ~~\forall ~i=0,1,\dots, M-1\},
\end{equation*}
where $\mathcal{P}^k$ denotes the polynomials of degree less or equal to $k\geq 1$. Moreover, let
\begin{equation*}
u(x_i^+)=\lim_{
\substack{\epsilon\rightarrow0 \\
\epsilon>0}}u(x_i+\epsilon), \qquad
u(x_i^-)=\lim_{
\substack{\epsilon\rightarrow0 \\
\epsilon>0}}u(x_i-\epsilon).
\end{equation*}
Then we define the jump and average of $u$ at the endpoints of $I_i$ as follows
  \begin{equation*}
\llbracket{u(x_i)}\rrbracket = u(x_i^-)-u(x_i^+),\quad
\boldsymbol{ \{}u(x_i)\boldsymbol{\}}=\frac{1}{2}\left(u(x_i^-)+u(x_i^+)\right)\quad \forall~ i=1,2,\dots, M-1.
 \end{equation*}
By convention, we also extend the definition of jump and average at the endpoints of the
unit interval
\begin{equation}\begin{aligned}\label{eq3.1}
&\llbracket{u(x_0)}\rrbracket = -u(x_0^+),\quad
~~\boldsymbol{ \{}u(x_0)\boldsymbol{\}}= u(x_0^+), \\
&\llbracket{u(x_M)}\rrbracket = u(x_M^-),\quad
~~\boldsymbol{ \{}u(x_M)\boldsymbol{\}}=u(x_M^-).
\end{aligned}\end{equation}

Multiplying \eqref{eq1.1} by $\xi(x)$ and integrating by parts on each interval $I_i$, with $\xi(x)$ being a function in $V_h$ for $i\in\{1,2,\dots, M-1\}$, we get
\begin{equation*}
\int_{x_i}^{x_{i+1}} \Delta u(x)\xi(x)dx=-\int_{x_i}^{x_{i+1}} \nabla u(x)\nabla \xi(x)dx
+\nabla u(x_{i+1}^-)\xi(x_{i+1}^-)-\nabla u(x_{i}^+)\xi(x_{i}^+).
\end{equation*}
Summing above equation from $i=0$ to $M-1$, and using \eqref{eq3.1}, we have
\begin{equation}\label{eq3.2}
\sum_{i=0}^{M-1}\int_{x_i}^{x_{i+1}} \Delta u(x)\xi(x)dx=-\sum_{i=0}^{M-1}\int_{x_i}^{x_{i+1}} \nabla u(x)\nabla \xi(x)dx + \sum_{i=0}^{M}\llbracket{\nabla u(x_i)\xi(x_i)}\rrbracket.
\end{equation}
It can be verified that 
\begin{equation*}
\llbracket{\nabla u(x_i)\xi(x_i)}\rrbracket=\boldsymbol{ \{}\nabla u(x_i)\boldsymbol{\}}\llbracket{\xi(x_i)}\rrbracket +\boldsymbol{ \{}\xi(x_i)\boldsymbol{\}}\llbracket{\nabla u(x_i)}\rrbracket, \quad 1\leq i\leq M-1.
\end{equation*}
Due to the fact that the exact solution satisfies $\boldsymbol{ \{}\xi(x_i)\boldsymbol{\}}\llbracket{\nabla u(x_i)}\rrbracket=0$ for $0\leq i\leq M$, \eqref{eq3.2} becomes
\begin{equation*}
\sum_{i=0}^{M-1}\int_{x_i}^{x_{i+1}} \Delta u(x)\xi(x)dx=-\sum_{i=0}^{M-1}\int_{x_i}^{x_{i+1}} \nabla u(x)\nabla \xi(x)dx + \sum_{i=0}^{M}\boldsymbol{ \{}\nabla u(x_i)\boldsymbol{\}}\llbracket{\xi(x_i)}\rrbracket.
\end{equation*}

Let $B_h$ be the symmetric interior penalty discrete bilinear form
\begin{equation}\begin{aligned}\label{eq3.3}
B_h(u, \xi)&=-\sum_{i=0}^{M-1}\int_{x_i}^{x_{i+1}} \nabla u(x)\nabla \xi(x)dx + \sum_{i=0}^{M}\boldsymbol{ \{}\nabla u(x_i)\boldsymbol{\}}\llbracket{\xi(x_i)}\rrbracket\\
&\quad+\sum_{i=0}^{M}\boldsymbol{ \{}\nabla \xi(x_i)\boldsymbol{\}}\llbracket{ u(x_i)}\rrbracket
-\sum_{i=0}^{M}\frac{\sigma}{h}\llbracket{u(x_i)}\rrbracket\llbracket{\xi(x_i)}\rrbracket,
\end{aligned}\end{equation}
where the positive constant $\sigma$ is the interior penalty stabilization parameter. 
 The third term in the right hand of \eqref{eq3.3} makes the bilinear form symmetric and the last term ensures coercivity of the bilinear with sufficiently large $\sigma$.
The bilinear form $B_h(\cdot, \cdot)$ defines a discrete linear operator 
$\Delta_h: V_h\rightarrow V_h$ as
\begin{equation*}
\langle\Delta_h U, \xi \rangle =B_h(U, \xi) \quad \forall~  \xi\in V_h, 
\end{equation*}
where $\langle\cdot,\cdot\rangle$ is the inner product in $L^2([a,b])$.

As a consequence, the discrete formulation is given by: find $ U,   V \in V_h$ such that
\begin{equation} \left\{\begin{aligned}\label{eq3.4}
d U &= Vdt,\\
d V&= \Delta_h U dt -P_hf( U)dt + P_hg(U)dW(t),\\
\end{aligned} \right. \end{equation} 
where $P_h: L^2([a,b])\rightarrow V_h$ is the projection.
Analogous to \eqref{eq2.2}, the semi-discrete scheme \eqref{eq3.4} preserves the following discrete version of energy evolution law.
\begin{thm} \label{thmDG}
The averaged energy for semi-discrete scheme \eqref{eq3.4} has the following evolutionary relationship
\begin{equation}\label{eq3.5}
\mathbb E\left[\widehat{H}(U(t),V(t))\right]=\mathbb E\left[\widehat{H}(U(0),V(0))\right]
+\frac{1}{2}\mathbb E\left[\int_0^tTr\left(P_hg(U)Q(P_hg(U))^*\right)ds\right],
\end{equation}
where
\begin{equation*}
\widehat{H}(U(t), V(t))=\frac{1}{2}\langle V(t), V(t) \rangle-\frac{1}{2}\langle \Delta_hU(t), U(t) \rangle+\int_{a}^{b}P_h \widetilde f(U(t))dx
\end{equation*}
is the discrete energy with respect to \eqref{eq1.2}.
\end{thm}
\noindent{\textbf{Proof. }}
By the  It\^o's  formula for $\frac{1}{2}\int_{a}^{b}V(t)^\top V(t)dx=\frac{1}{2}\langle V(t), V(t) \rangle$, we deduce
\begin{equation}\begin{aligned}\label{eq3.6}
\frac{1}{2}\langle V(t), V(t) \rangle&=\frac{1}{2}\langle V(0), V(0) \rangle
+ \int_{0}^{t}\langle V(s), P_hg(U(s)) dW(s)\rangle+\int_{0}^{t}\langle V(s), \Delta_h U(s) -P_hf(U(s)) \rangle ds\\
&\quad+\frac{1}{2}\int_{0}^{t} Tr\big[(P_hg(U)Q^{\frac{1}{2}})(P_hg(U)Q^{\frac{1}{2}})^*\big]ds.
\end{aligned}\end{equation}
Using the first equation of \eqref{eq3.4}, we get
\begin{equation}\begin{aligned}\label{eq3.7}
&\quad \int_{0}^{t}\langle V(t), \Delta_h U(s) -P_hf(U(s)) \rangle ds\\
&=\int_{0}^{t}\left \langle \frac{d U(s)}{ds}, \Delta_h U(s)\right \rangle ds-\int_{a}^{b}\int_{0}^{t}\frac{d U(s)}{ds}P_hf(U(s))dsdx\\
&=\frac{1}{2}\langle \Delta_h U(s),  U(s) \rangle\Big|_0^t-\int_{a}^{b}P_h\widetilde f(U(s))\Big|_0^tdx\\
&=\frac{1}{2}\langle \Delta_h U(t),  U(t) \rangle-\frac{1}{2}\langle \Delta_hU(0), U(0) \rangle-\int_{a}^{b}P_h\widetilde f(U(t))dx+\int_{a}^{b}P_h\widetilde f( U(0))dx.
\end{aligned}\end{equation}
We combine \eqref{eq3.6}, \eqref{eq3.7} and take expectation to obtain 
\begin{equation}\begin{aligned}\label{eq3.8}
&\quad \mathbb E\left[\frac{1}{2}\langle V(t), V(t) \rangle-\frac{1}{2}\langle \Delta_hU(t), U(t) \rangle+\int_{a}^{b}P_h\widetilde f(U(t))dx\right]\\
&=\mathbb E\left[\frac{1}{2}\langle V(0), V(0) \rangle-\frac{1}{2}\langle \Delta_h U(0), U(0) \rangle+\int_{a}^{b}P_h\widetilde f( U(0))dx\right]\\
&\quad+\frac{1}{2}\mathbb E\left[\int_0^tTr\left(P_hg(U)Q(P_hg(U))^*\right)ds\right].
\end{aligned}\end{equation}
Then according to the definition of $\widehat{H}(U(t), V(t))$, we derive
\begin{equation*}
\mathbb E\left[\widehat{H}(U(t),V(t))\right]=\mathbb E\left[\widehat{H}(U(0),V(0))\right]
+\frac{1}{2}\mathbb E\left[\int_0^tTr\left(P_hg(U)Q(P_hg(U))^*\right)ds\right],
\end{equation*}
which completes the proof.
\hfill $\square$\\

\begin{rem}
 If we consider \eqref{eq1.1} with periodic boundary condition,  the corresponding jump and average at the endpoints of the unit interval are
\begin{equation*}\begin{aligned}
&\llbracket{u(x_0)}\rrbracket_p = u(x_M^-)-u(x_0^+),\quad
~~\boldsymbol{ \{}u(x_0)\boldsymbol{\}}_p= \frac{1}{2}\left(u(x_0^+)+u(x_M^-)\right), \\
&\llbracket{u(x_M)}\rrbracket_p = u(x_M^-)-u(x_0^+),\quad
\boldsymbol{ \{}u(x_M)\boldsymbol{\}}_p=\frac{1}{2}\left(u(x_0^+)+u(x_M^-)\right).
\end{aligned}\end{equation*}
In addition, the symmetric interior penalty discrete bilinear form is
\begin{equation*}\begin{aligned}
B_{hp}(u, \xi)&=\sum_{i=0}^{M-1}\int_{x_i}^{x_{i+1}} \Delta u(x)\xi(x)dx=-\sum_{i=0}^{M-1}\int_{x_i}^{x_{i+1}} \nabla u(x)\nabla \xi(x)dx\\
&\quad + \sum_{i=0}^{M-1}\boldsymbol{ \{}\nabla u(x_i)\boldsymbol{\}}_p\llbracket{\xi(x_i)}\rrbracket_p
+\sum_{i=0}^{M-1}\boldsymbol{ \{}\nabla \xi(x_i)\boldsymbol{\}}_p\llbracket{ u(x_i)}\rrbracket_p
+\sum_{i=0}^{M-1}\frac{\sigma}{h}\llbracket{u(x_i)}\rrbracket_p\llbracket{\xi(x_i)}\rrbracket_p.
\end{aligned}\end{equation*}
Moreover, the corresponding semi-discrete scheme discretized by the interior penalty discontinuous Galerkin finite element method is to find $ U,   V \in \overline{V}_{h}$ such that 
\begin{equation} \left\{\begin{aligned}\label{Peq2}
d U &= Vdt,\\
d V&= \Delta_{hp} U dt -P_{hp}f( U)dt + P_{hp}g(U)dW(t),\\
\end{aligned} \right. \end{equation} 
where $\overline{V}_{h}=\{u\in L^2([a, b]): u|_{I_i} \in \mathcal{P}^k,~u(a)=u(b) ~~\forall ~i=0,1,\dots, M-1\}$,
$P_{hp}: L^2([a,b])\rightarrow \overline{V}_{h}$ is the projection, and $\Delta_{hp}: \overline{V}_{h}\rightarrow \overline{V}_{h}$ is a discrete linear operator defined by bilinear form $B_h(\cdot, \cdot)$ as follows
\begin{equation*}
\langle\Delta_{hp} U, \xi \rangle =B_{hp}(U, \xi) \quad \forall~  \xi\in \overline{V}_{h}. 
\end{equation*}
Similarly,  the semi-discrete scheme \eqref{Peq2}  of the nonlinear stochastic wave equation with periodic boundary condition preserves the averaged energy evolution law, i.e.,
\begin{equation*}
\mathbb E\left[\check{H}(U(t),V(t))\right]=\mathbb E\left[\check{H}(U(0),V(0))\right]
+\frac{1}{2}\mathbb E\left[\int_0^tTr\left(P_{hp}g(U)Q(P_{hp}g(U))^*\right)ds\right],
\end{equation*}
where
\begin{equation*}
\check{H}(U(t), V(t))=\frac{1}{2}\langle V(t), V(t) \rangle-\frac{1}{2}\langle \Delta_{hp}U(t), U(t) \rangle+\int_{a}^{b}P_{hp} \widetilde f(U(t))dx.
\end{equation*}
\end{rem}

%----------------------全离散--------------------------------------------------------------------
\section{Energy-preserving  fully-discrete schemes}
In this section, we turn to considering  energy-preserving fully-discrete schemes for the nonlinear stochastic wave equation  \eqref{eq1.1} based on  semi-discrete schemes in Section 2.  What follows  mainly focuses on the construction of  fully-discrete schemes by discretizing \eqref{eq2.2} in temporal direction. We denote by the $n$th time level $t_n=n\Delta t$, $n=0,1,\dots, N$  with time step size $\Delta t = T/N$.

We first rewrite the semi-discrete scheme \eqref{eq2.2} as 
\begin{align*}
d\mbf X(t)= \widetilde{\mbf A}\mbf X(t)dt+\mbf J\mbf \Phi (\mbf X(t))dt+\mbf \Upsilon(\mbf X(t)) d\boldsymbol {\beta}(t),
\end{align*}
where 
\begin{equation*}
\mbf X=
\begin{bmatrix}
\mbf U \\ \mbf V
\end{bmatrix},
\quad 
\widetilde{\mbf A}=\begin{bmatrix}
0 & \mbf I\\
\mbf{A^{-1}D} & 0
\end{bmatrix},\quad 
\mbf J=\begin{bmatrix}
0 & \mbf I\\
-\mbf I & 0
\end{bmatrix},\quad 
\mbf \Phi (\mbf X) =\begin{bmatrix}
\mbf F(\mbf U)\\
0
\end{bmatrix},\quad
\mbf \Upsilon (\mbf X)=\begin{bmatrix}
0\\
G(\mbf U)\mbf e \boldsymbol{\Lambda}
\end{bmatrix}.
\end{equation*}
It can be known that the solution given by the constant variation method reads
\begin{equation}
\label{mildsolu}
\mbf X(t)=\exp(\widetilde{\mbf A}t)\mbf X(0)+\int_0^t \exp(\widetilde{\mbf A}(t-s))\mbf J
\mbf \Phi (\mbf X(s))ds+\int_0^t \exp(\widetilde{\mbf A}(t-s))\mbf \Upsilon( \mbf X(s))d\boldsymbol {\beta}(s).
\end{equation}

 The Pad\'e approximation of a rational function is given by ratio of two polynomials. The coefficients of the polynomial in both the numerator and the denominator are determined by using the coefficients in the Taylor series expansion of the function. The main advantage of the Pad\'e approximation over the Taylor series approximation is that the Taylor series approximation can exhibit oscillation which may produce an approximation error bound (see \cite{Ongun}).  In what follows, we use the Pad\'e approximation to approximate the exponential matrix to construct the fully-discrete schemes which can preserve the averaged energy evolution law.

As is well known, the matrix exponential $\exp{(\mbf R)}$ for $(M-1)$-dimensional matrix $\mbf R$ has the Taylor expansion
\begin{equation*}
\exp{(\mbf R)}=\mbf I+\sum_{i=1}^{\infty} \frac{\mbf R^{i}}{i!}.
\end{equation*}
A simple way to approximate the exponential function $\exp{(x)}$ is making use of the rational Pad\'e approximation
\begin{equation}
\exp{(x)} \approx P_{(r,s)}(x)= D_{(r,s)}^{-1}(x)N_{(r,s)}(x),
\end{equation}
where
$$D_{(r,s)}(x)=1+\sum_{i=1}^{s} \frac{(r+s-i)!s!}{(r+s)!i!(s-i)!}(-x)^{i}=1+\sum_{i=1}^{s}\tilde{a}_i(-x)^{i},$$
$$N_{(r,s)}(x)=1+\sum_{i=1}^{r} \frac{(r+s-i)!r!}{(r+s)!i!(r-i)!}x^{i}=1+\sum_{i=1}^{r}a_ix^{i}$$
with 
\begin{equation*}
\tilde{a}_i=\frac{(r+s-i)!s!}{(r+s)!i!(s-i)!}, \qquad a_i=\frac{(r+s-i)!r!}{(r+s)!i!(r-i)!}.
\end{equation*}

By exploiting the Pad\'e approximation of the exponential matrix and combining with the discrete gradient method, we propose a methodology of constructing the energy-preserving numerical scheme as follows
\begin{equation}\label{scheme 2}
\mbf X^{n+1}=P_{l,l}(\Delta t \widetilde{\mbf A})\mbf X^n+\widetilde{\mbf A}^{-1}(P_{l,l}(\Delta t \widetilde{\mbf A})-I)\mbf J \mbf \Phi (\mbf X^{n},\mbf X^{n+1})+P_{l,l}(\Delta t \widetilde{\mbf A})\mbf \Upsilon(\mbf X^n) \Delta\boldsymbol {\beta}_n
\end{equation}
with 
\begin{equation*}
P_{l,l}(\Delta t \widetilde{\mbf A})=D_{l,l}^{-1}(\Delta t \widetilde{\mbf A})N_{l,l}(\Delta t \widetilde{\mbf A}),
\quad \mbf \Phi (\mbf X^{n},\mbf X^{n+1}) =\begin{bmatrix}
\overline{\nabla}_U F(\mbf U^{n+1},\mbf U^{n})\\
0
\end{bmatrix},\quad
\mbf \Upsilon (\mbf X^n)=\begin{bmatrix}
0\\
G(\mbf U^n)\mbf e \boldsymbol{\Lambda}
\end{bmatrix},
\end{equation*}
where
\begin{equation*}\begin{aligned}
\overline{\nabla}_U F(\mbf U^{n+1},\mbf U^{n})&=\left(\frac{\widetilde f(U_1^{n+1})-\widetilde f(U_1^{n})}{U_1^{n+1}-U_1^{n}}, \frac{\widetilde f(U_2^{n+1})-\widetilde f(U_2^{n})}{U_2^{n+1}-U_2^{n}}, \dots,\frac{\widetilde f(U_{M-1}^{n+1})-\widetilde f(U_{M-1}^{n})}{U_{M-1}^{n+1}-U_{M-1}^{n}}\right)^\top,\\
\Delta \boldsymbol {\beta}_n&=\left(\beta_1(t_{n+1})-\beta_1(t_{n}), \beta_2(t_{n+1})-\beta_2(t_{n}), \dots,\beta_P(t_{n+1})-\beta_P(t_{n})\right)^\top .
\end{aligned} \end{equation*}
 In order to deal with the unboundedness of $\Delta_j \beta_n=\beta_j(t_{n+1})-\beta_j(t_{n})$ for $j\in \{1,2,\dots,P\}$, which is simulated by $\sqrt{\Delta t}   \xi_j^n$ with $\xi_j^n  \sim \mathcal N(0,1)$, we introduce the following truncated random variable $\Delta_j \widehat{\beta}_n=\widehat{\beta}_j(t_{n+1})-\widehat{\beta}_j(t_{n}) = \sqrt{\Delta t}   \widehat{\xi}_j^n$ 
\begin{equation*} \widehat{\xi}_j^n = \left\{\begin{aligned}
&\quad \xi_j^n,\quad\qquad if~ |\xi_j^n|\leq A_{\Delta t},\\
&\quad A_{\Delta t},~\qquad if ~\xi_j^n > A_{\Delta t},\\
& -A_{\Delta t},\quad\quad if ~\xi_j^n < -A_{\Delta t}, \\
\end{aligned} \right. \end{equation*} 
with $A_{\Delta t}:=\sqrt{2k|ln\Delta t|}$, where $k\geq 2$ is an integer. Then  $|\widehat{\xi}_j^n |\leq \sqrt{2k|ln\Delta t|}$. To avoid confusion, the truncated random variable $\Delta_j \widehat{\beta}_n$ is still denoted by $\Delta_j \beta_n$.

More specifically,  we can rewrite \eqref{scheme 2} as 
\begin{equation} \label{eq2.12}\left\{\begin{aligned}
\mbf U^{n+1}-\mbf U^{n}&= -\sum_{k=1}^{[l/2]}  b_k(\mbf U^{n+1}-\mbf U^{n})+\sum_{k=0}^{[(l-1)/2]}  c_k(\mbf V^{n+1}+\mbf V^{n})\\
 &+\sum_{k=0}^{[(l-1)/2]}  c_kG(\mbf U^n)\mbf E \boldsymbol{\Lambda} \Delta\boldsymbol {\beta}_n,\\
\mbf V^{n+1}-\mbf V^{n} &= -\sum_{k=1}^{[l/2]}  b_k(\mbf V^{n+1}-\mbf V^{n})+\sum_{k=0}^{[(l-1)/2]}  c_k\mbf{A^{-1}D}(\mbf U^{n+1}+\mbf U^{n})\\
&-2\sum_{k=0}^{[(l-1)/2]}  c_k\overline{\nabla}_U F(\mbf U^{n+1},\mbf U^{n})+\sum_{k=0}^{[l/2]}  b_kG(\mbf U^n)\mbf E \boldsymbol{\Lambda} \Delta\boldsymbol {\beta}_n,
\end{aligned} \right. \end{equation} 
where $b_k=a_{2k}(\Delta t)^{2k}(\mbf{A^{-1}D})^k$ and $c_k=a_{2k+1}(\Delta t)^{2k+1}(\mbf{A^{-1}D})^k$.

When we take $l=1,2$, \eqref{eq2.12} becomes
\begin{equation} \label{P11}\left\{\begin{aligned}
\mbf U^{n+1}-\mbf U^{n}&= \frac{\Delta t}{2}(\mbf V^{n+1}+\mbf V^{n})+\frac{\Delta t}{2}G(\mbf U^n)\mbf E \boldsymbol{\Lambda} \Delta\boldsymbol {\beta}_n,\\
\mbf V^{n+1}-\mbf V^{n} &= \frac{\Delta t}{2}\mbf{A^{-1}D}(\mbf U^{n+1}+\mbf U^{n})
-\Delta t\overline{\nabla}_U F(\mbf U^{n+1},\mbf U^{n})+G(\mbf U^n)\mbf E \boldsymbol{\Lambda} \Delta\boldsymbol {\beta}_n,
\end{aligned} \right. \end{equation} 
and
\begin{equation}\label{P22} \left\{\begin{aligned}
\mbf U^{n+1}-\mbf U^{n}&= \frac{\Delta t}{2}(\mbf V^{n+1}+\mbf V^{n})-\frac{\Delta t^2}{12}\mbf{A^{-1}D}(\mbf U^{n+1}-\mbf U^{n})+\frac{\Delta t}{2}G(\mbf U^n)\mbf E \boldsymbol{\Lambda} \Delta\boldsymbol {\beta}_n,\\
\mbf V^{n+1}-\mbf V^{n} &= \frac{\Delta t}{2}\mbf{A^{-1}D}(\mbf U^{n+1}+\mbf U^{n})
-\Delta t\overline{\nabla}_U F(\mbf U^{n+1},\mbf U^{n})-\frac{\Delta t^2}{12}\mbf{A^{-1}D}(\mbf V^{n+1}-\mbf V^{n})\\
&\quad+G(\mbf U^n)\mbf E \boldsymbol{\Lambda} \Delta\boldsymbol {\beta}_n +\frac{\Delta t^2}{12}\mbf{A^{-1}D}G(\mbf U^n)\mbf E \boldsymbol{\Lambda} \Delta\boldsymbol {\beta}_n.
\end{aligned} \right. \end{equation}  

It is obvious that the fully-discrete scheme \eqref{eq2.12} is implicit. The unique solvability can be proved similarly as Lemma 2.4 in \cite{Milstein1} by using the contraction mapping principle, the global Lipschitz condition and the fact that $-\mbf A^{-1}\mbf D$ is a symmetric positive definite matrix. We now prove that the fully-discrete scheme of the nonlinear stochastic wave equation, that
is, the numerical solution given by \eqref{eq2.12} satisfies the averaged energy evolution law.
\begin{thm} \label{thm4}
The averaged energy for fully-discrete scheme \eqref{eq2.12} has the following evolutionary relationship
\begin{equation}\label{eq2.13}
\mathbb E\left[\bar{H}(\mbf U^{n+1},\mbf V^{n+1})\right]=\mathbb E\left[\bar{H}(\mbf U^n,\mbf V^n)\right]
+\frac{h\Delta t}{2}\mathbb E\left[Tr( G^2(\mbf U^n)\mbf E \boldsymbol{\Lambda}^2\mbf E^\top) \right].
\end{equation}
\end{thm}	
\noindent{\textbf{Proof. }}
For the sake of simplicity, we assume that $l=2m+1,$ with $m\geq 1, $ since  the case that $l$ is an even integer is similar. Based on \eqref{eq2.12}, we deduce
%It can be verified that
%\begin{equation*}\begin{aligned}
%D_{l,l}(\Delta t \widetilde{\mbf A})&=\sum_{k=0}^m a_{2k}(\Delta t)^{2k}\widetilde{\mbf A}^{2k}-\sum_{k=0}^m a_{2k+1}(\Delta t)^{2k+1}\widetilde{\mbf A}^{2k+1},\\
%N_{l,l}(\Delta t \widetilde{\mbf A})&=\sum_{k=0}^m a_{2k}(\Delta t)^{2k}\widetilde{\mbf A}^{2k}+\sum_{k=0}^m a_{2k+1}(\Delta t)^{2k+1}\widetilde{\mbf A}^{2k+1}.\\
%\end{aligned}\end{equation*}
%Denoting $b_k=a_{2k}(\Delta t)^{2k}(\mbf{A^{-1}D})^k$ and $c_k=a_{2k+1}(\Delta t)^{2k+1}(\mbf{A^{-1}D})^k$,  we deduce 
\begin{equation}\label{eq2.14} \begin{array}{ll}
 \begin{bmatrix}
	\ds \mbf{U}^{n+1}-\mbf{U}^{n}\\[0.1in]
	\ds \mbf{V}^{n+1}-\mbf{V}^{n}
	\end{bmatrix} 
&= \begin{bmatrix}
	\ds -\sum_{k=1}^m b_k(\mbf{U}^{n+1}-\mbf{U}^{n})+\sum_{k=0}^m c_k(\mbf{V}^{n+1}+\mbf{V}^{n})\\[0.1in]
	\ds -\sum_{k=1}^m b_k(\mbf{V}^{n+1}-\mbf{V}^{n})+\sum_{k=0}^m c_k\mbf{A^{-1}D}(\mbf{U}^{n+1}+\mbf{U}^{n})-2\sum_{k=0}^mc_k\overline{\nabla}_U F(\mbf U^{n+1},\mbf U^{n})
	\end{bmatrix}\\[0.2in]
&\quad +\begin{bmatrix}
	\ds \sum_{k=0}^m c_kG(\mbf U^n)\mbf E \boldsymbol{\Lambda} \Delta \boldsymbol {\beta}_n\\[0.1in]
	\ds \sum_{k=0}^m b_kG(\mbf U^n)\mbf E \boldsymbol{\Lambda} \Delta \boldsymbol {\beta}_n
	\end{bmatrix}.\\[0.1in]
	\end{array}
\end{equation}
Multiplying 
$[(\mbf U^{n+1}+\mbf U^n)^\top, 
	(\mbf V^{n+1}+\mbf V^n)^\top]\begin{bmatrix}
	-\mbf{A^{-1}D}& 0\\
	0& \mbf I
	\end{bmatrix}$ on the both side of \eqref{eq2.14}, we have 
\begin{equation}\begin{aligned}\label{eq2.15}
&\quad -\sum_{i=1}^{M-1}\left(U_i^{n+1}\mathcal A^{-1} \delta_x^2U_i^{n+1}
-U_i^n\mathcal A^{-1} \delta_x^2U_i^n\right)+\sum_{i=1}^{M-1}\left((V_i^{n+1})^2 -(V_i^{n} )^2\right)\\
&=(\mbf U^{n+1}+\mbf U^n)^\top\mbf{A^{-1}D}\sum_{k=1}^m b_k(\mbf U^{n+1}-\mbf U^n)\\
&\quad+(\mbf V^{n+1}+\mbf V^n)^\top\Big[-\sum_{k=1}^m b_k(\mbf V^{n+1}-\mbf V^n)-2\sum_{k=0}^mc_k\overline{\nabla}_U F(\mbf U^{n+1},\mbf U^{n})\Big]\\
&
\quad-(\mbf U^{n+1}+\mbf U^n)^\top\mbf{A^{-1}D}\sum_{k=0}^m c_kG(\mbf U^n) \mbf E \boldsymbol{\Lambda} \Delta \boldsymbol {\beta}_n+
(\mbf V^{n+1}+\mbf V^n)^\top\sum_{k=0}^m b_kG(\mbf U^n) \mbf E \boldsymbol{\Lambda} \Delta \boldsymbol {\beta}_n\\
&:=I_1+I_2+I_3+I_4,
\end{aligned}\end{equation}
where $U_i^n$ and $V_i^n$ are the $i$th component of vectors $\mbf U^n$ and $\mbf V^n$ for $i\in\{1,2,\dots, M-1\}$.
Based on \eqref{eq2.14}, we derive 
\begin{align}\label{eq2.16}
I_1+I_2
&=(\mbf U^{n+1}+\mbf U^n)^\top\mbf{A^{-1}D}\sum^m_{k=1} b_k(\sum_{k=0}^m  b_k)^{-1}(\sum\limits_{k=0}^m c_k)(\mbf V^{n+1}+\mbf V^n)\nonumber\\
&\quad-(\mbf V^{n+1}+\mbf V^n)^\top\sum^m_{k=1} b_k(\sum_{k=0}^m  b_k)^{-1}(\sum_{k=0}^m c_k)\mbf{A^{-1}D}(\mbf U^{n+1}+\mbf U^n)\nonumber\\
&\quad+2(\mbf V^{n+1}+\mbf V^n)^\top\sum^m_{k=1} b_k(\sum_{k=0}^m  b_k)^{-1}(\sum_{k=0}^m c_k)\overline{\nabla}_U F(\mbf U^{n+1},\mbf U^{n})\nonumber\\
&\quad- 2(\mbf V^{n+1}+\mbf V^n)^\top(\sum_{k=0}^m c_k)\overline{\nabla}_U F(\mbf U^{n+1},\mbf U^{n})\\
&\quad+(\mbf U^{n+1}+\mbf U^n)^\top\mbf{A^{-1}D}\sum_{k=1}^mb_k(\sum_{k=0}^m b_k)^{-1}(\sum_{k=0}^m c_k)G(\mbf U^n) \mbf E \boldsymbol{\Lambda} \Delta \boldsymbol {\beta}_n\nonumber\\
&\quad-(\mbf V^{n+1}+\mbf V^n)^\top\sum^m_{k=1} b_k G(\mbf U^n) \mbf E \boldsymbol{\Lambda} \Delta \boldsymbol {\beta}_n\nonumber\\
&=: II_1+II_2+II_3+II_4+II_5+II_6.\nonumber
\end{align}
The symmetry of $\mbf{A^{-1}D}$ leads to
\begin{equation}\label{eq2.17}
II_1+II_2=0,
\end{equation}
and
\begin{equation}
II_3+II_4 = -2(\mbf V^{n+1}+\mbf V^n)^\top(\sum_{k=0}^m  b_k)^{-1}(\sum_{k=0}^m c_k)\overline{\nabla}_U F(\mbf U^{n+1},\mbf U^{n}).
\end{equation}
The first equation of \eqref{eq2.14} equals to 
\begin{align*}
\mbf U^{n+1}-\mbf U^n=(\sum_{k=0}^m  b_k)^{-1}(\sum_{k=0}^m c_k)(\mbf V^{n+1}+\mbf V^n)
+(\sum_{k=0}^m  b_k)^{-1}(\sum_{k=0}^m c_k)G(\mbf U^n) \mbf E \boldsymbol{\Lambda} \Delta \boldsymbol {\beta}_n,
\end{align*}
which yields 
\begin{align*}
II_3+II_4=-2\sum_{i=1}^{M-1}\left(\widetilde f(U_i^{n+1})-\widetilde f(U_i^{n})\right)+2(G(\mbf U^n) \mbf E \boldsymbol{\Lambda} \Delta \boldsymbol {\beta}_n)^\top(\sum_{k=0}^m  b_k)^{-1}(\sum_{k=0}^m c_k)\overline{\nabla}_U F(\mbf U^{n+1},\mbf U^{n}).
\end{align*}
As a consequence,
\begin{equation}\begin{aligned}\label{eq2.18}
&\quad I_3+I_4+II_3+II_4+II_5+II_6\\
&=-2\sum_{i=1}^{M-1}\left(\widetilde f(U_i^{n+1})-\widetilde f(U_i^{n})\right)+2(G(\mbf U^n) \mbf E \boldsymbol{\Lambda} \Delta \boldsymbol {\beta}_n)^\top(\sum_{k=0}^m  b_k)^{-1}(\sum_{k=0}^m c_k)\overline{\nabla}_U F(\mbf U^{n+1},\mbf U^{n})
\\
&\quad-(\mbf U^{n+1}+\mbf U^n)^\top\mbf{A^{-1}D}\sum_{k=0}^m c_kG(\mbf U^n) \mbf E \boldsymbol{\Lambda}\Delta \boldsymbol {\beta}_n+
(\mbf V^{n+1}+\mbf V^n)^\top\sum_{k=0}^m b_kG(\mbf U^n) \mbf E\boldsymbol{\Lambda} \Delta \boldsymbol {\beta}_n\\
&\quad+(\mbf U^{n+1}+\mbf U^n)^\top\mbf{A^{-1}D}\sum_{k=1}^mb_k(\sum_{k=0}^m b_k)^{-1}(\sum_{k=0}^m c_k)G(\mbf U^n) \mbf E \boldsymbol{\Lambda} \Delta \boldsymbol {\beta}_n\\
&\quad-(\mbf V^{n+1}+\mbf V^n)^\top\sum^m_{k=1} b_k G(\mbf U^n) \mbf E \boldsymbol{\Lambda} \Delta \boldsymbol {\beta}_n\\
&=-2\sum_{i=1}^{M-1}\left(\widetilde f(U_i^{n+1})-\widetilde f(U_i^{n})\right)+2(G(\mbf U^n) \mbf E \boldsymbol{\Lambda} \Delta \boldsymbol {\beta}_n)^\top(\sum_{k=0}^m  b_k)^{-1}(\sum_{k=0}^m c_k)\overline{\nabla}_U F(\mbf U^{n+1},\mbf U^{n})\\
&\quad-(\mbf U^{n+1}+\mbf U^n)^\top\mbf{A^{-1}D}(\sum_{k=0}^m b_k)^{-1}(\sum_{k=0}^m c_k)G(\mbf U^n) \mbf E \boldsymbol{\Lambda}\Delta \boldsymbol {\beta}_n\\
&\quad+(\mbf V^{n+1}+\mbf V^n)^\top G(\mbf U^n) \mbf E \boldsymbol{\Lambda}\Delta  \boldsymbol {\beta}_n.
\end{aligned}\end{equation}
By using the second equation of  \eqref{eq2.14},
\begin{equation*}\begin{aligned}
(\mbf V^{n+1}-\mbf V^n)&=(\sum_{k=0}^m b_k)^{-1}(\sum_{k=0}^m c_k)\mbf{A^{-1}D}(\mbf U^{n+1}+\mbf U^n)\\
&\quad-2(\sum_{k=0}^m b_k)^{-1}(\sum_{k=0}^m c_k)\overline{\nabla}_U F(\mbf U^{n+1},\mbf U^{n})+G(\mbf U^n) \mbf E \boldsymbol{\Lambda} \Delta\boldsymbol {\beta}_n.
\end{aligned}\end{equation*}
We have the estimate of \eqref{eq2.18} as follows
\begin{align}\label{eq2.19}
\nonumber
&\quad I_3+I_4+II_3+II_4+II_5+II_6\\
\nonumber
&=-2\sum_{i=1}^{M-1}\left(\widetilde f(U_i^{n+1})-\widetilde f(U_i^{n})\right)-(G(\mbf U^n) \mbf E \boldsymbol{\Lambda} \Delta\boldsymbol {\beta}_n)^\top(\mbf V^{n+1}-\mbf V^n)\\
\nonumber
&\quad+(G(\mbf U^n) \mbf E \boldsymbol{\Lambda}\Delta\boldsymbol {\beta}_n)^\top(\sum_{k=0}^m b_k)^{-1}(\sum_{k=0}^m c_k)\mbf{A^{-1}D}(\mbf U^{n+1}+\mbf U^n)\\
&\quad+(G(\mbf U^n) \mbf E \boldsymbol{\Lambda} \Delta\boldsymbol {\beta}_n)^\top G(\mbf U^n)\mbf E \boldsymbol{\Lambda} \Delta\boldsymbol {\beta}_n\\
\nonumber
&\quad-(\mbf U^{n+1}+\mbf U^n)^\tau\mbf{A^{-1}D}(\sum_{k=0}^m b_k)^{-1}(\sum_{k=0}^m c_k)G(\mbf U^n) \mbf E \boldsymbol{\Lambda} \Delta\boldsymbol {\beta}_n\\
\nonumber
&\quad+(\mbf V^{n+1}+\mbf V^n)^\top G(\mbf U^n) \mbf E \boldsymbol{\Lambda} \Delta\boldsymbol {\beta}_n.
\end{align}

We combine the above estimates of \eqref{eq2.15}-\eqref{eq2.17} with \eqref{eq2.19} and take expectation to obtain  
\begin{equation}\begin{aligned}\label{eq2.20}
&\quad \mathbb E\left[-\sum_{i=1}^{M-1}\left(U_i^{n+1}\mathcal A^{-1} \delta_x^2U_i^{n+1}
-U_i^n\mathcal A^{-1} \delta_x^2U_i^n\right)+\sum_{i=1}^{M-1}\left((V_i^{n+1})^2 -(V_i^{n} )^2\right)\right]\\
&=\mathbb E\left[-2\sum_{i=1}^{M-1}\left(\widetilde f(U_i^{n+1})-\widetilde f(U_i^{n})\right)\right) + \mathbb E\left((G(\mbf U^n) \mbf E \boldsymbol{\Lambda}\Delta \boldsymbol {\beta}_n)^\top G(\mbf U^n)\mbf E \boldsymbol{\Lambda} \Delta \boldsymbol {\beta}_n\right].
\end{aligned}\end{equation}
Multiplying 
$\frac{h}{2}$ on the both side of \eqref{eq2.20} and noting the definition of $\bar{H}(\mbf U, \mbf V)$, we have
\begin{equation*}
\mathbb E\left[\bar{H}(\mbf U^{n+1}, \mbf V^{n+1})\right]=\mathbb E\left[\bar{H}(\mbf U^{n}, \mbf V^{n})\right]+\frac{h\Delta t}{2}\mathbb E\left[Tr( \mbf E G^2(\mbf U^n)\boldsymbol{\Lambda}^2\mbf E^\top \right],
\end{equation*}
which finishes the proof.
\hfill $\square$\\

\begin{rem}
Similarly, for the semi-discrete scheme \eqref{eq3.4}  which is discretized by the interior penalty discontinuous Galerkin finite element method for solving  \eqref{eq1.1} with Dirichlet boundary condition, we have the fully-discrete scheme
\begin{equation} \label{eq2.22}\left\{\begin{aligned}
U^{n+1}-U^{n}&= -\sum_{k=1}^{[l/2]}  \tilde{b}_k(U^{n+1}-U^{n})+\sum_{k=0}^{[(l-1)/2]}  \tilde{c}_k(V^{n+1}+V^{n})\\
 &\quad +\sum_{k=0}^{[(l-1)/2]}  \tilde{c}_kP_hg(U^n)\Delta W_n,\\
V^{n+1}-V^{n} &= -\sum_{k=1}^{[l/2]}  \tilde{b}_k(V^{n+1}-V^{n})+\sum_{k=0}^{[(l-1)/2]}  \tilde{c}_k\Delta_h(U^{n+1}+U^{n})\\
&\quad -2\sum_{k=0}^{[(l-1)/2]}  \tilde{c}_kP_h\frac{\widetilde f(U^{n+1})-\widetilde f(U^{n})}{U^{n+1}-U^{n}}+\sum_{k=0}^{[l/2]}  \tilde{b}_kP_hg(U^n)\Delta W_n,
\end{aligned} \right. \end{equation} 
where $\tilde{b}_k=a_{2k}(\Delta t)^{2k}\Delta_h^k$ and $\tilde{c}_k=a_{2k+1}(\Delta t)^{2k+1}\Delta_h^k$, and $\Delta W_n= W(t_{n+1})- W(t_n)$.

 In addition, using the same analysis technique in Theorem \ref{thm4}, the averaged energy for fully-discrete scheme \eqref{eq2.22} has the following relationship
\begin{equation*}
\mathbb E\left[\hat{H}(U^{n+1},V^{n+1})\right]=\mathbb E\left[\hat{H}(U^n,V^n)\right]
+\frac{\Delta t}{2}\mathbb E\left[\int_a^b Tr(P_hg(U^n)Q (P_hg(U^n))^*)dx\right].
\end{equation*}
\end{rem}

\begin{rem}
For the semi-discrete schemes \eqref{Peq1} and  \eqref{Peq2}, which approximate the nonlinear stochastic wave equation \eqref{eq1.1} with periodic boundary condition, two fully-discrete schemes \eqref{eq2.12}  and \eqref{eq2.22}  can be given corresponding, respectively. Besides, the unique solvability and  averaged energy evolution law of fully-discrete schemes can be given by using similar arguments.
\end{rem}
%---------------------------------------------------------------------------------------------------------------
\begin{figure}
	\centering
	\subfigure{
		\begin{minipage}{6cm}
			\centering
			\includegraphics[height=5cm,width=6cm]{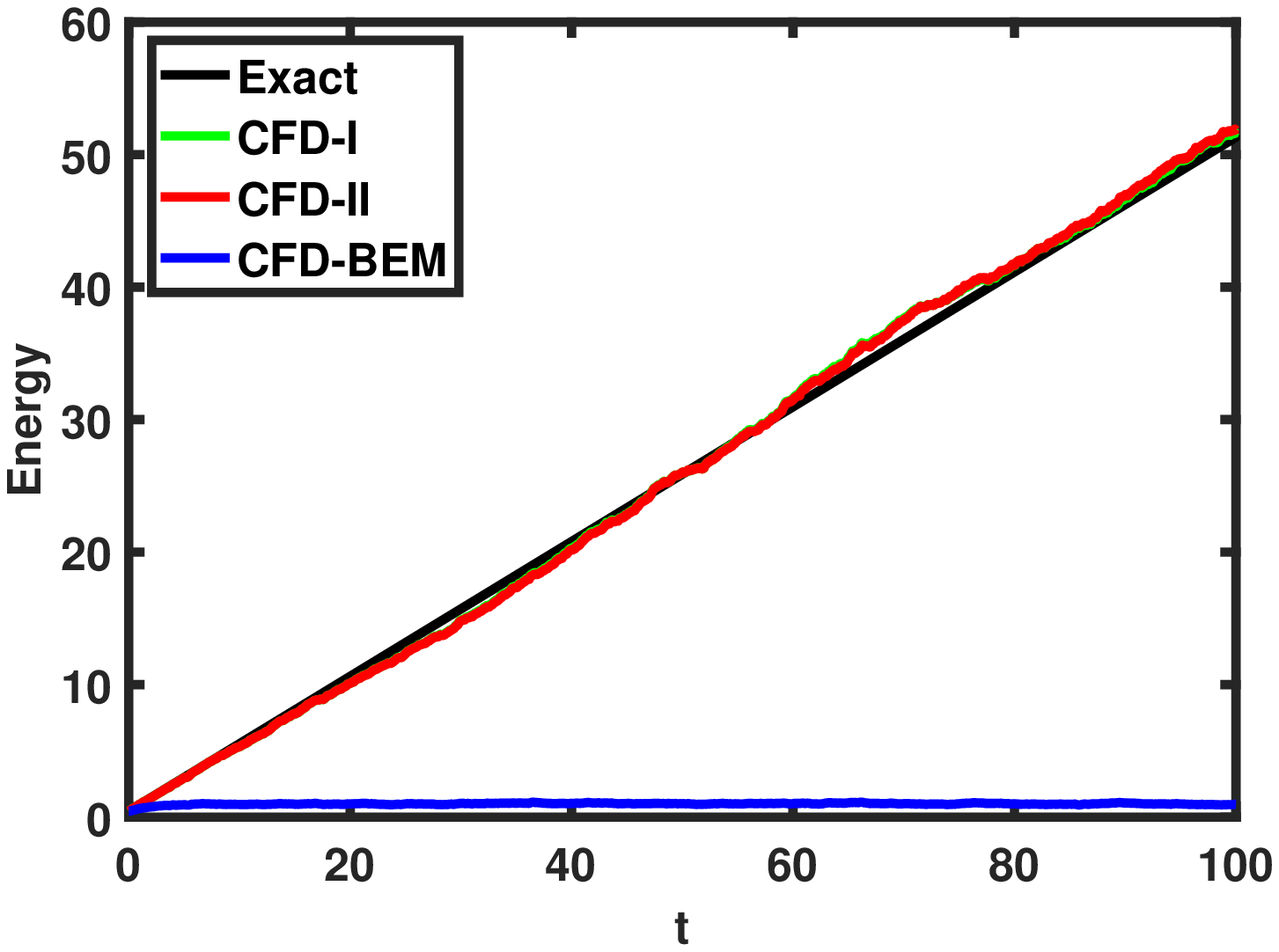}
		\end{minipage}
	}
	\subfigure{
		\begin{minipage}{6cm}
			\centering
			\includegraphics[height=5cm,width=6cm]{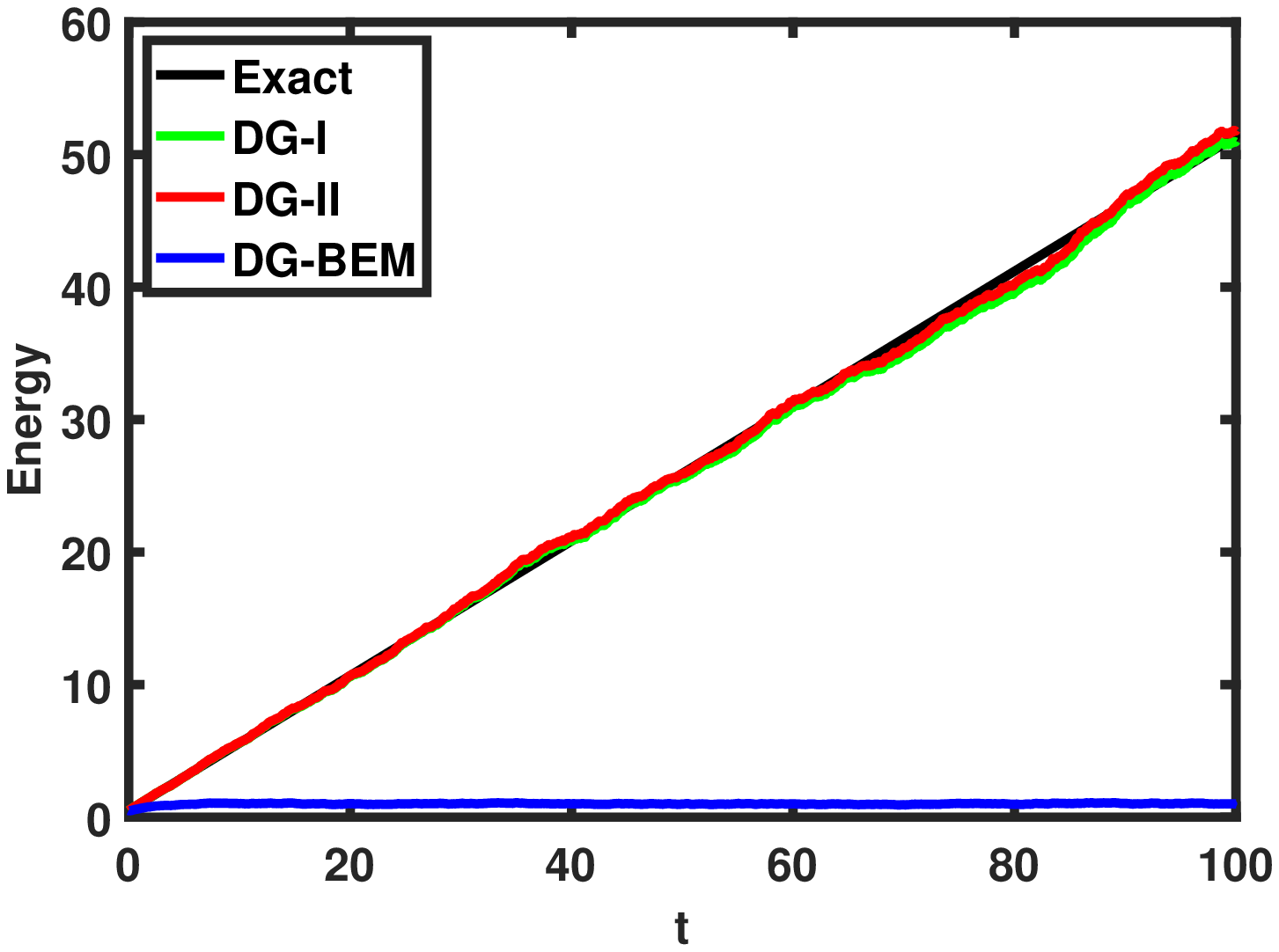}
		\end{minipage}
	}
	\caption{ Averaged energy evolution relationship ($f(u) = 0, g(u) = 1$) with $\Delta t= 1/20, h=1/10$}
	\label{fig1}
\end{figure}

%------------数值算例---------------------------------
\section{Numerical experiments}

This section presents various numerical experiments in order to illustrate the energy-preserving property of the proposed fully-discrete  schemes \eqref{eq2.12} and \eqref{eq2.22} with $l = 1, 2$ for the 1-dimensional nonlinear stochastic wave equation under the homogeneous Dirichlet boundary condition.  
We will compare the proposed numerical schemes with the following classical fully-discrete schemes based on\\

{\noindent 1. backward} Euler-Maruyama (BEM) method (\cite{Kloeden}) and compact finite difference method \eqref{eq2.2}
\begin{equation}\label{eq9.1}\left\{\begin{aligned}
\mbf U ^{n+1}&= \mbf U ^{n} + \Delta t \mbf V^{n+1},\\
\mbf V ^{n+1}&= \mbf V ^{n} + \Delta t  \mbf{A^{-1}D}\mbf U ^{n+1} -\Delta t  F(\mbf U ^{n+1})dt 
+ G(\mbf U^n)\mbf e \boldsymbol{\Lambda}\Delta \boldsymbol {\beta}_n,\\
\end{aligned} \right. \end{equation} 
2. backward Euler-Maruyama (BEM)  method (\cite{Kloeden}) and discontinuous Galerkin finite element method \eqref{eq3.4}
\begin{equation}\label{eq9.3} \left\{\begin{aligned}
U^{n+1} &= U^{n}+\Delta t V^{n+1},\\
V^{n+1} &= V^n+\Delta t\Delta_hU^{n+1}-\Delta tP_hf( U^{n+1})+P_hg(U^n)\Delta W_n,
\end{aligned} \right. \end{equation} 
3. semi-implicit Crank-Nicolson-Maruyama (CNM)  method (\cite{Hausenblas2, Walsh2}) and compact finite difference method \eqref{eq2.2}
\begin{equation}\label{eq9.2}\left\{\begin{aligned}
\mbf U ^{n+1}&= \mbf U ^{n} + \Delta t \frac{\mbf V^n+\mbf V^{n+1}}{2},\\
\mbf V ^{n+1}&= \mbf V ^{n} + \Delta t\mbf{A^{-1}D}\frac{\mbf U ^{n}+\mbf U ^{n+1}}{2} -\Delta t  F(\mbf U ^n) 
+ G(\mbf U^n)\mbf e \boldsymbol{\Lambda}\Delta \boldsymbol {\beta}_n,\\
\end{aligned} \right. \end{equation} 	
4. semi-implicit Crank-Nicolson-Maruyama (CNM) method (\cite{Hausenblas2, Walsh2}) and discontinuous Galerkin finite element method \eqref{eq3.4}
\begin{equation} \label{eq9.4}\left\{\begin{aligned}
U^{n+1} &= U^{n}+\Delta t\frac{V^n + V^{n+1}}{2},\\
V^{n+1} &= V^n+\Delta t\Delta_h\frac{U^{n+1}+U^{n}}{2}-\Delta tP_hf( U^n)+P_hg(U^n)\Delta W_n.
\end{aligned} \right. \end{equation} 
For convenience, we denote numerical schemes \eqref{P11}, \eqref{P22}, \eqref{eq9.1} and \eqref{eq9.2} by CFD-I, CFD-II, CFD-BEM and CFD-CNM, respectively, in the case that the semi-discrete scheme is based on the compact finite difference method. 
Similarly, for the case that discontinuous Galerkin finite element method is applied in spatial direction, the numerical schemes \eqref{eq2.22} with $l=1, 2$, \eqref{eq9.3}, \eqref{eq9.4} are denoted by DG-I, DG-II, DG-BEM and DG-CNM, respectively.  
The first and second examples (Subections \ref{sec;ne1} and \ref{sec;ne2}) test the numerical approximation by simulating the
stochastic wave equation with globally Lipschitz continuous coefficients. In Subection \ref{sec;ne3}, numerical tests of the proposed fully-discrete schemes
for the non-globally Lipschitz case are presented. In all the experiments, the expectation is approximated
by taking average over 1000 realizations.

\begin{figure}
	\centering
	\subfigure{
		\begin{minipage}{6cm}
			\centering
			\includegraphics[height=5cm,width=6cm]{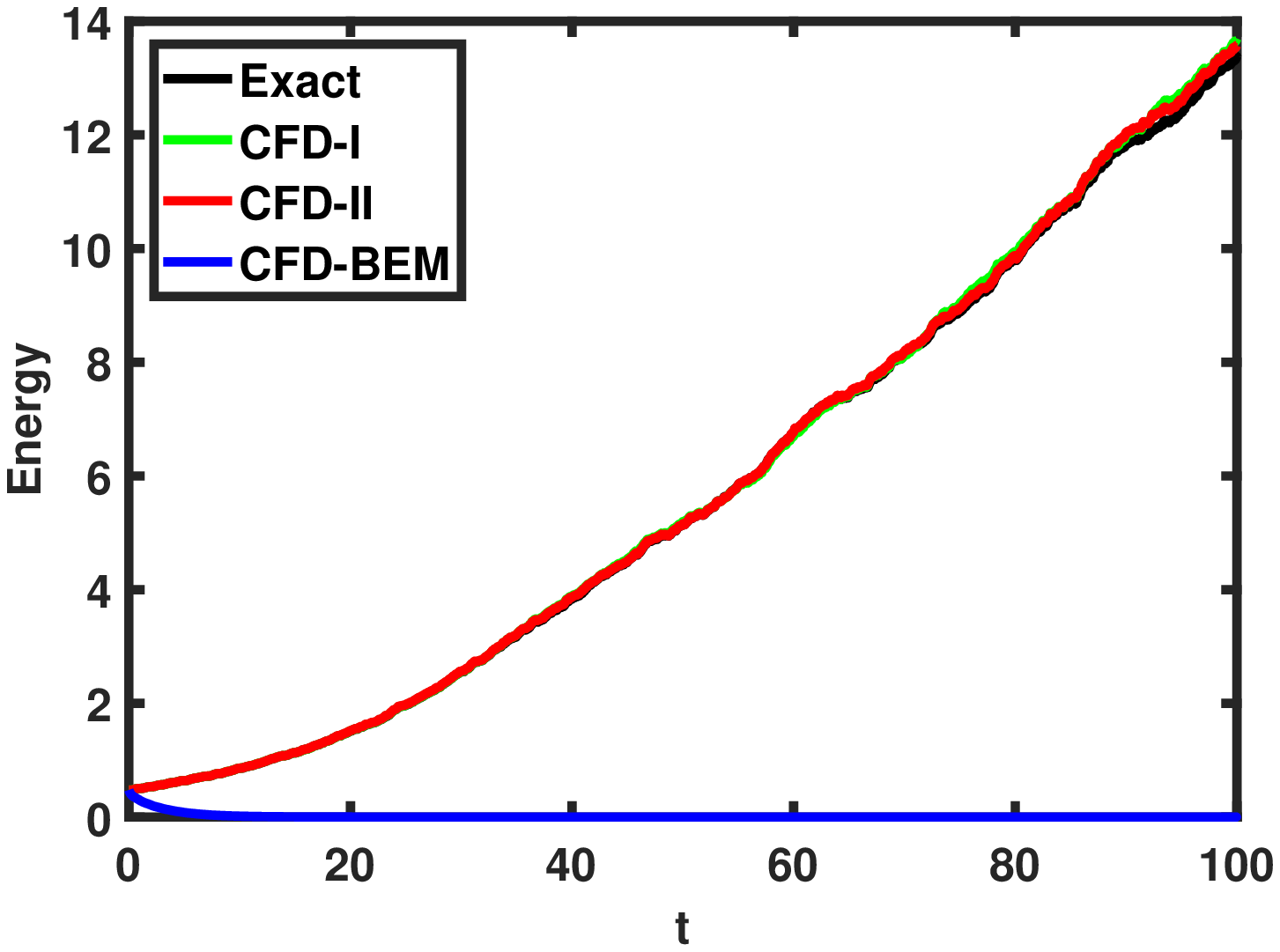}
		\end{minipage}
	}
	\subfigure{
		\begin{minipage}{6cm}
			\centering
			\includegraphics[height=5cm,width=6cm]{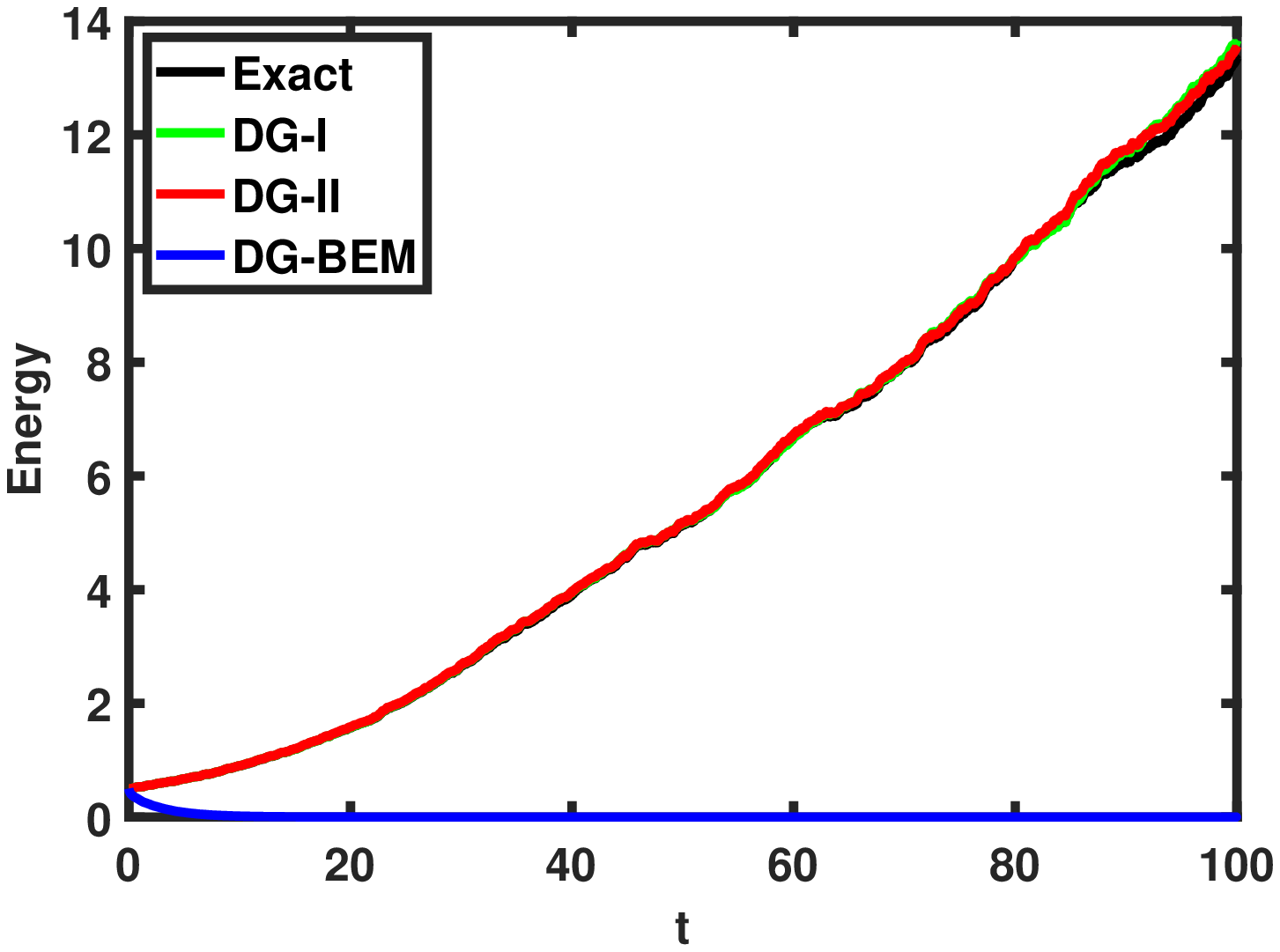}
		\end{minipage}
	}
	\caption{ Averaged energy evolution relationship ($f(u) = 0, g(u) = \sin (u)$) with $\Delta t= 1/25, h=1/20$}
	\label{fig2}
\end{figure}

\begin{figure}
	\centering
	\subfigure{
		\begin{minipage}{6cm}
			\centering
			\includegraphics[height=5cm,width=6cm]{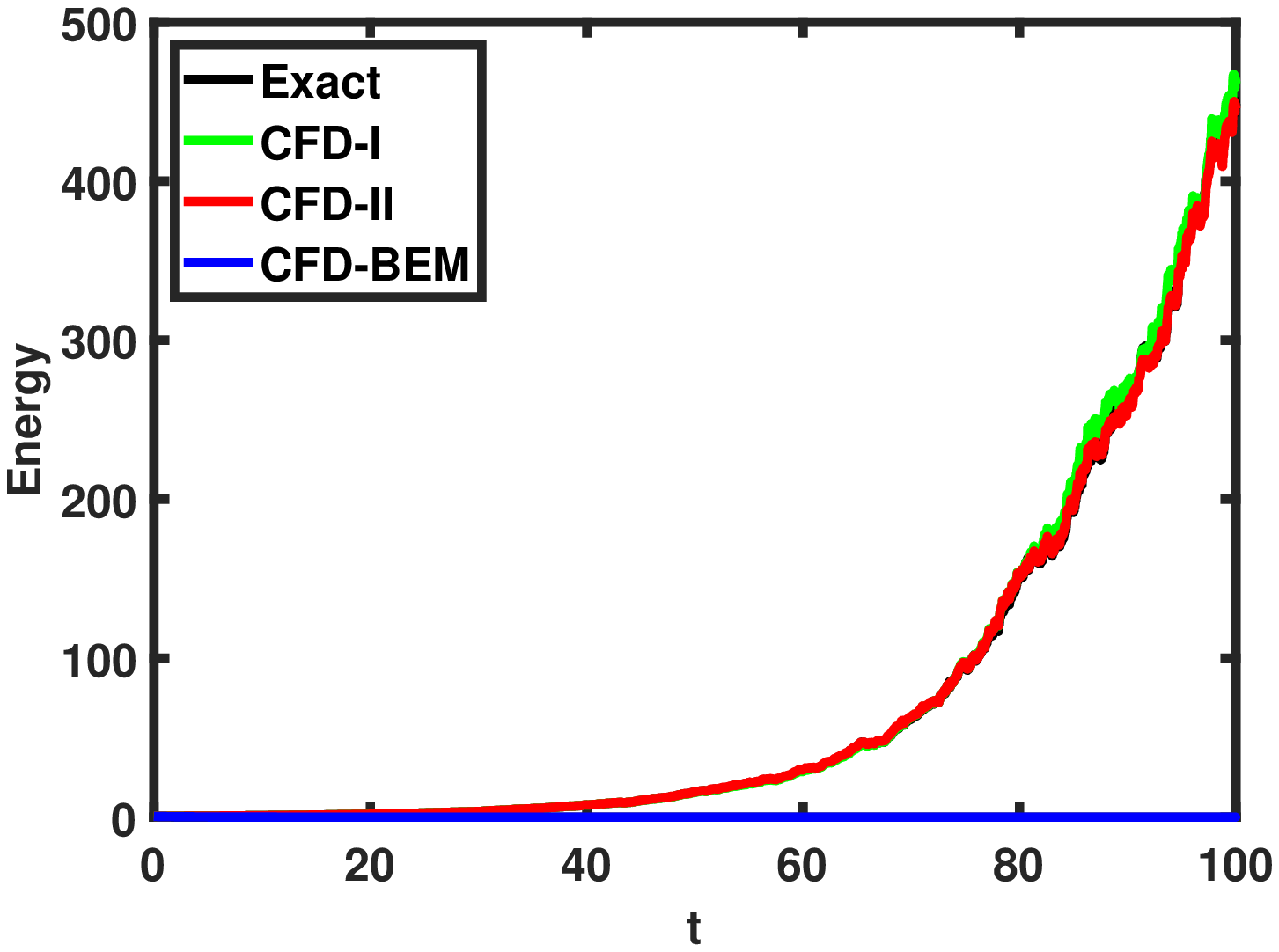}
		\end{minipage}
	}
	\subfigure{
		\begin{minipage}{6cm}
			\centering
			\includegraphics[height=5cm,width=6cm]{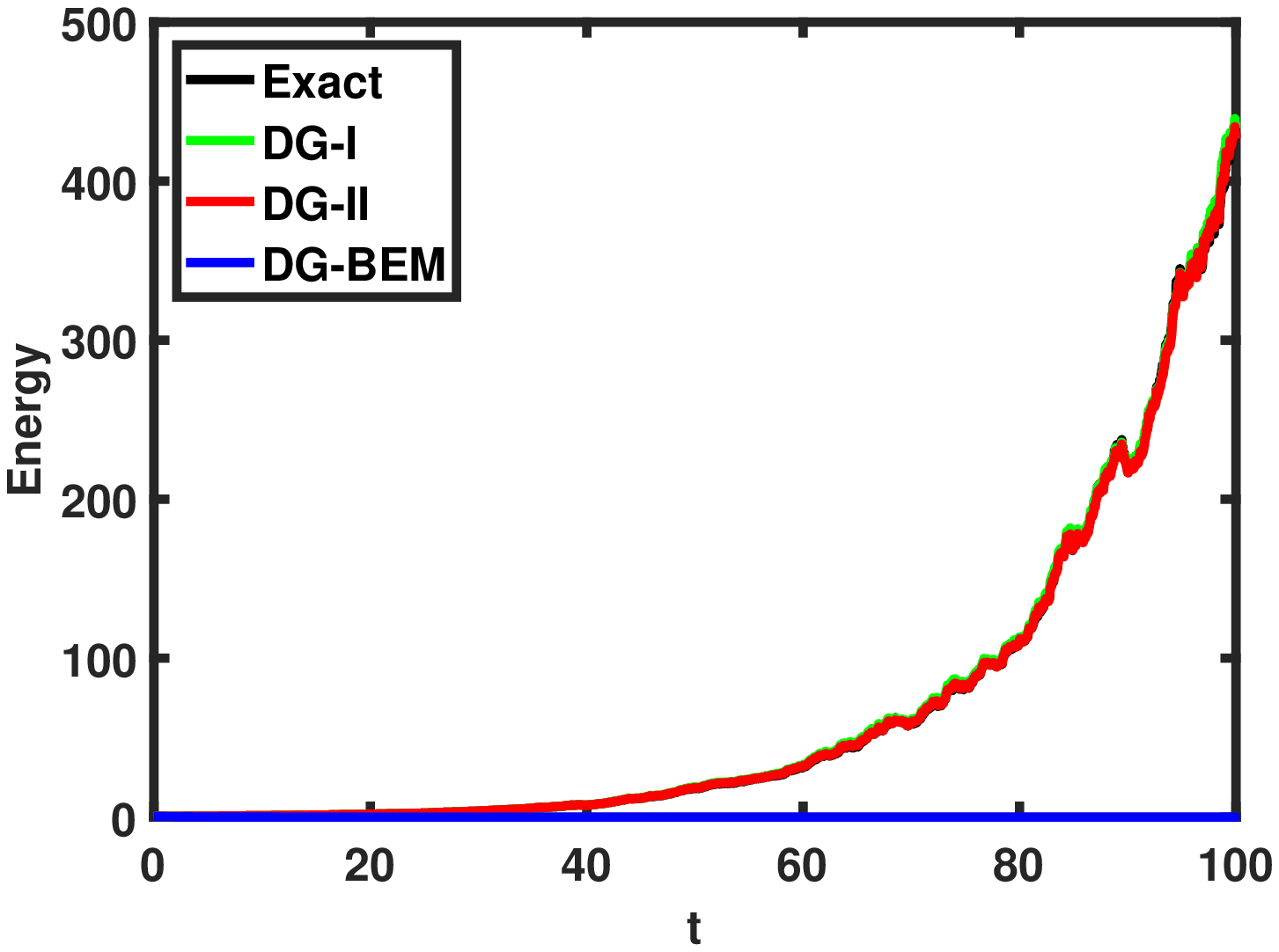}
		\end{minipage}
	}
	\caption{ Averaged energy evolution relationship ($f(u) = 0, g(u) = u$) with $\Delta t= 1/30, h=1/20$}
	\label{fig3}
\end{figure}
%-----------------------------------------------------
\subsection{Stochastic wave equation with $f(u)=0$ } 
\label{sec;ne1}
 Consider the following stochastic wave equation
\begin{equation} \label{eq6.1}\left\{\begin{aligned}
&du = vdt,  \\
&dv = \Delta u dt  +g(u) dW(t), &&(x, t) \in (0, 1) \times (0,100]
\end{aligned} \right. \end{equation}
 with initial conditions $u(x,0)=0, v(x,0)=1$. 
In this experiment, the diffusion coefficients are chosen as  $g(u)=1$, $g(u)=\sin(u)$ and $g(u)=u$,  which correspond to the additive noise case and multiplicative noise cases, respectively.
In the sequel, we choose the orthonormal basis $\left\{e_{k}\right\}_{k\in \mathbb{N}+}$ and the corresponding eigenvalue $\left\{q_{k}\right\}_{k\in \mathbb{N}+}$ of $Q$ as
\begin{align*}
e_{k}=\sqrt{2} \sin (k \pi x), \quad q_{k}=\frac{1}{k^{6}}.
\end{align*}

When $g(u) =1,$ it is known that the averaged energy of the exact solution grows linearly as time increases. 
The left hand picture of Fig. \ref{fig1} plots the quantity $\frac{h}{2}\mathbb E[(\mbf V^n)^\top \mbf V^n- (\mbf U^n)^\top \mbf{A^{-1}D}\mbf U^n]$ for CFD-I, CFD-II and CFD-BEM, and the right hand picture of Fig. \ref{fig1} plots the quantity $\frac{1}{2}\mathbb E[\langle V^n, V^n \rangle-\langle \Delta_h U^n, U^n \rangle]$ for DG-I, DG-II and DG-BEM for $n =1,\ldots,N,$ in the case of the stochastic wave equation \eqref{eq6.1} with additive noise. 
Moreover, the reference straight line (black line) in Fig. \ref{fig1} stands for the averaged energy evolution law of the exact solution, and has slope $\frac{1}{2}\sum_{k=1}^{N}q_k$ with $N = 2000$.
It can be observed that the proposed schemes named CFD-I, CFD-II, DG-I and DG-II reproduce the linear growth of the averaged energy, but both CFD-BEM and DG-BEM do not. 
The numerical results coincide with the fact that the BEM method applied does not preserve the averaged evolution law of the linear stochastic oscillator (see \cite{Strommen}).  
Fig. \ref{fig2} and Fig. \ref{fig3} demonstrate the evolution of the discrete averaged energy for CFD-I, CFD-II, CFD-BEM, DG-I, DG-II and DG-BEM of the stochastic wave equation \eqref{eq6.1} with $g(u) = \sin(u)$ and $g(u) = u$, respectively. 
We can see that the averaged energy in Fig. \ref{fig3} grows faster than that in Fig. \ref{fig2}.
Furthermore, it can be also checked that the proposed four schemes preserve perfectly the averaged energy evolution law, while CFD-BEM and DG-BEM fail.

%---------------------------------------------------------------------------------------------------------------
\begin{figure}
	\centering
	\subfigure{
		\begin{minipage}{6cm}
			\centering
			\includegraphics[height=5cm,width=6cm]{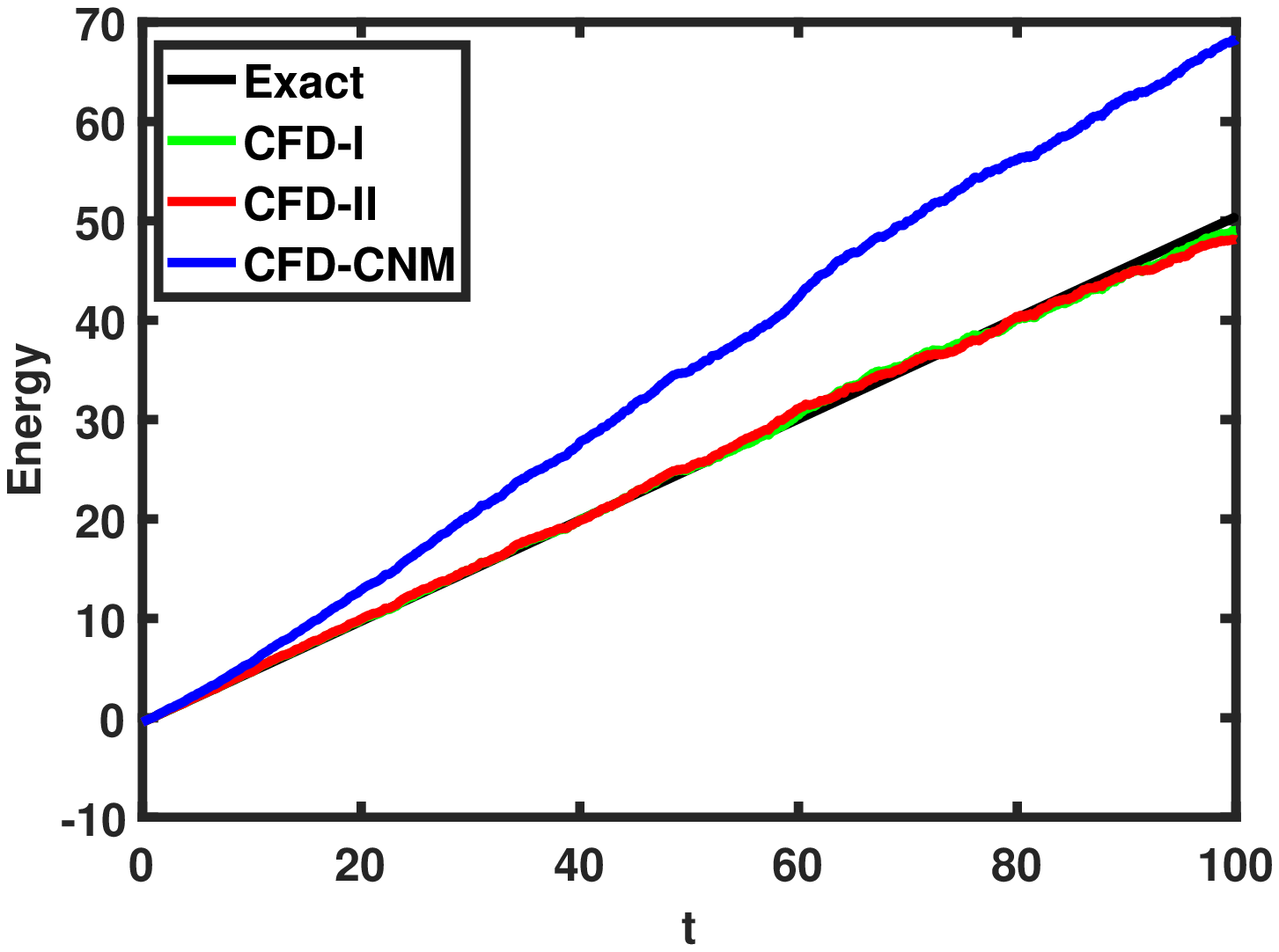}
		\end{minipage}
	}
	\subfigure{
		\begin{minipage}{6cm}
			\centering
			\includegraphics[height=5cm,width=6cm]{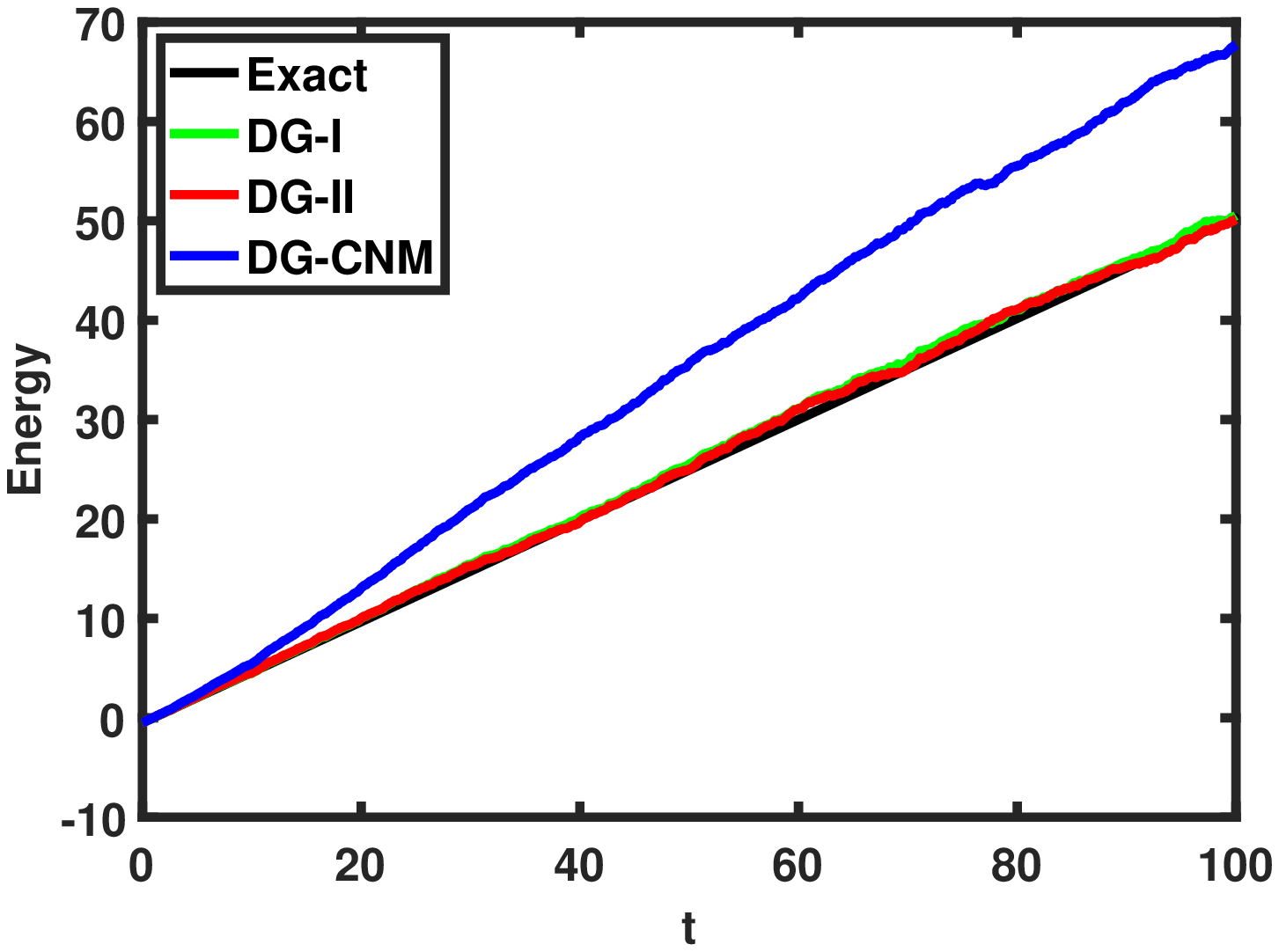}
		\end{minipage}
	}
	\caption{ Averaged energy evolution relationship ($f(u) = \sin(u), g(u) = 1$) with $\Delta t= 1/12, h=1/10$}
	\label{fig4}
\end{figure}
\subsection{Nonlinear stochastic wave equation with $f(u)=\sin (u)$  } 
\label{sec;ne2}
 In this section, we consider the following  nonlinear stochastic equation
\begin{equation} \label{eq6.2}\left\{\begin{aligned}
&du = vdt,  \\
&dv = \Delta u dt - \sin(u) dt +g(u) dW(t), &&(x, t) \in (0, 1) \times (0,100]
\end{aligned} \right. \end{equation}
with $u(x,0)=0, v(x,0)=1$, and take the cases $g(u)=1$ and $g(u)=\sin (u)$ into account. 

\begin{figure}
\centering
\subfigure{
\begin{minipage}{6cm}
\centering
\includegraphics[height=5cm,width=6cm]{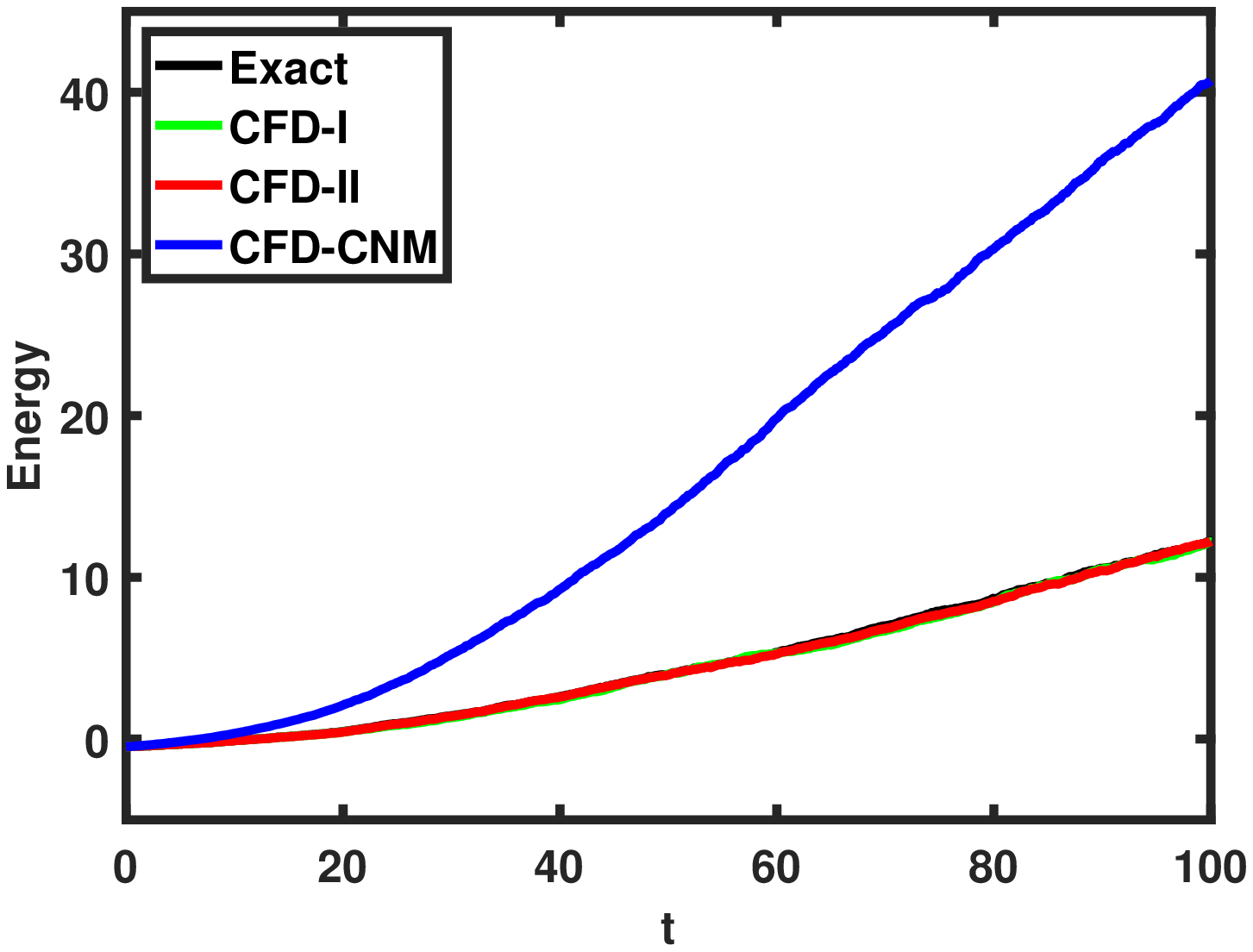}
\end{minipage}
}
\subfigure{
\begin{minipage}{6cm}
\centering
\includegraphics[height=5cm,width=6cm]{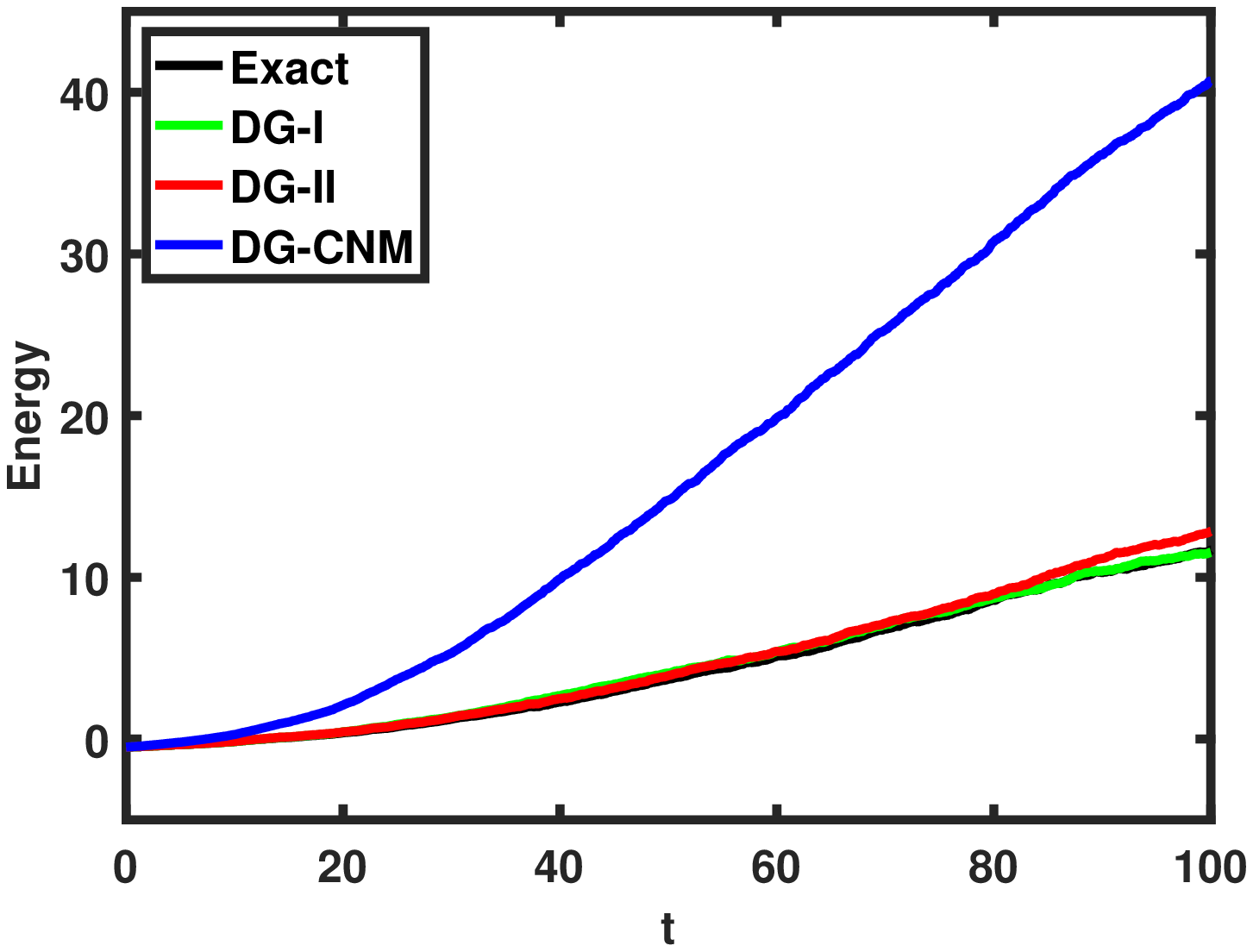}
\end{minipage}
}
\caption{Averaged energy evolution relationship ($f(u) = \sin(u), g(u) = \sin (u)$) with $\Delta t= 1/10, h=1/20$}
\label{fig5}
\end{figure}

\begin{figure}
\centering
\subfigure{
\begin{minipage}{6cm}
\centering
\includegraphics[height=5cm,width=6cm]{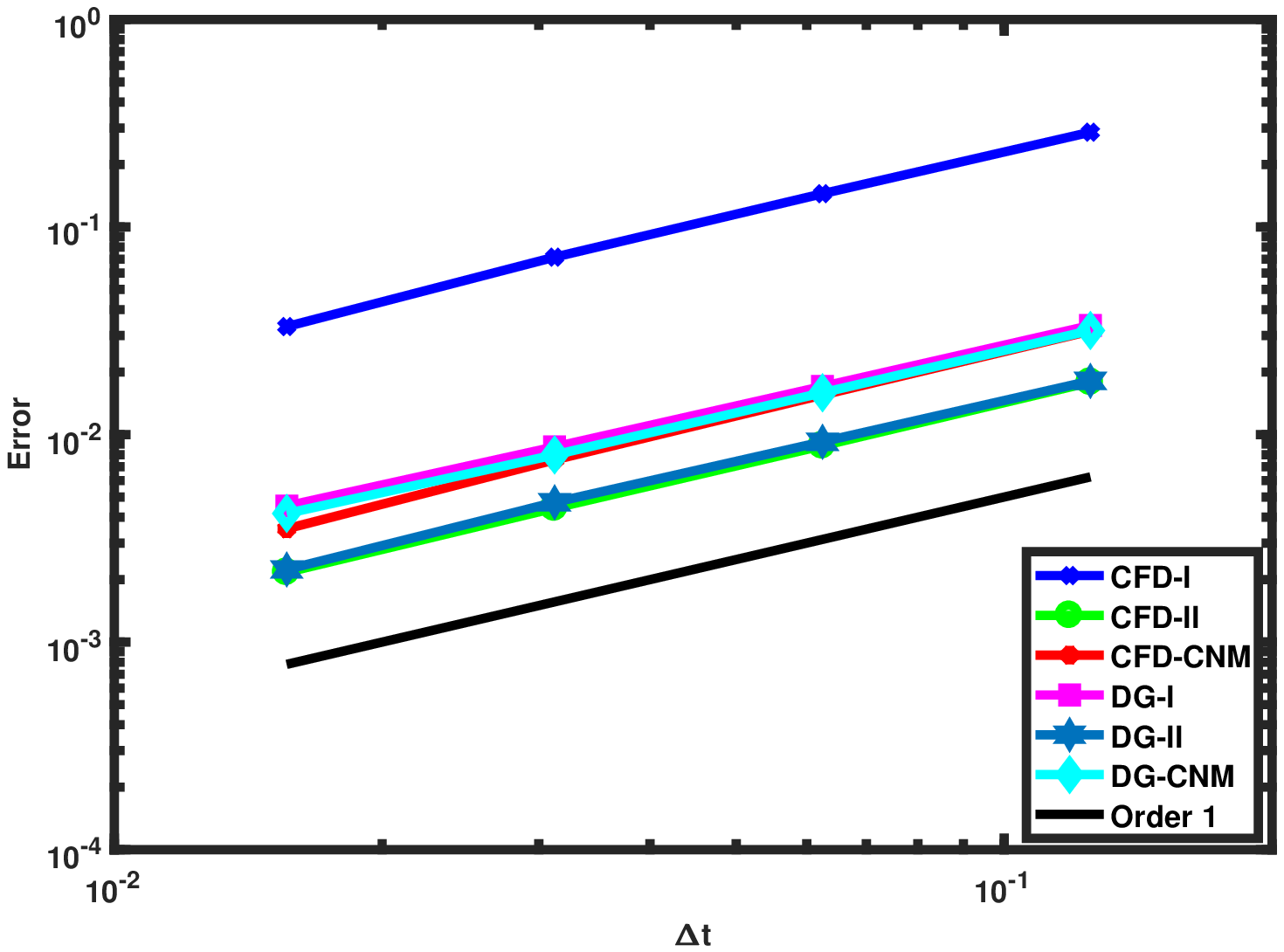}
\end{minipage}
}
\subfigure{
\begin{minipage}{6cm}
\centering
\includegraphics[height=5cm,width=6cm]{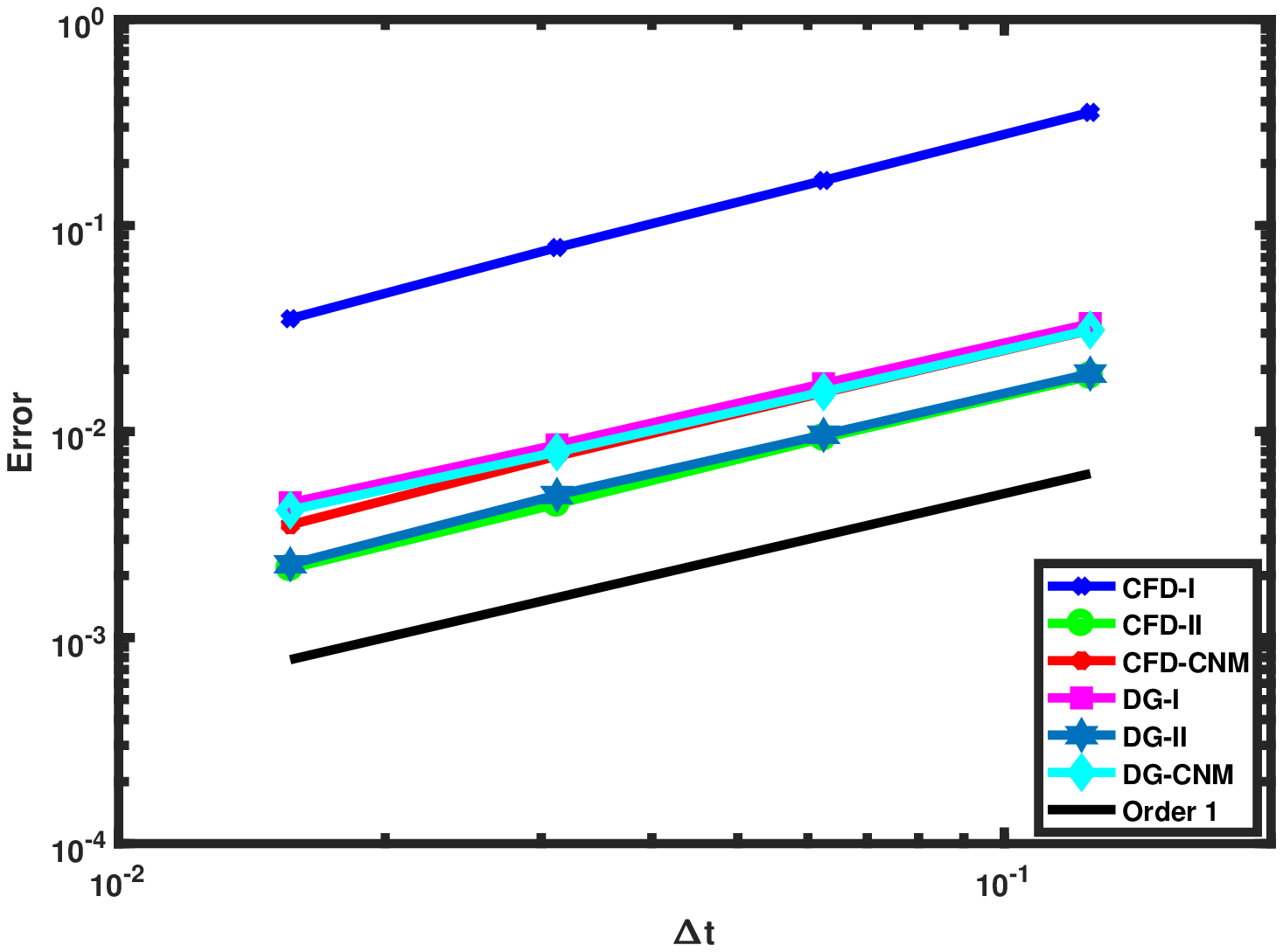}
\end{minipage}
}
\caption{Temporal rates of convergence: (left) $f(u) = \sin(u), g(u) = \sin (u)$, (right)  $f(u) = u^3, g(u) = \sin (u)$ }
\label{fig6.5}
\end{figure}

Figs. \ref{fig4}-\ref{fig5} present the evolution of discrete averaged energies for fully-discrete schemes named CFD-I, CFD-II, DG-I, DG-II, CFD-CNM and DG-CNM. 
From Fig. \ref{fig4}, it can be seen that  averaged energies associated with numerical solutions of  CFD-I, CFD-II, DG-I and DG-II grow linearly with the time raising, and coincide  with the averaged energy of exact solution when the stochastic wave equation is driven by additive noise. 
Although the expected energies of CFD-CNM and DG-CNM also possess the linear growth property,  the slope of blue lines is greater than the one of the reference black one, which means that  CFD-CNM and DG-CNM could not inherit the averaged energy evolution law of the original nonlinear stochastic wave equation with additive noise. 
In the case of stochastic wave equation \eqref{eq6.1} driven by multiplicative noise, from Fig. \ref{fig5} it also can be found the preservation of averaged energy of the proposed schemes, and the unpreservation of CFD-CNM and DG-CNM. 
The numerical results are consistent with the theoretical results. 
Fig. \ref{fig6.5} displays the temporal approximation errors $\|u(T)-U^N\|$ against $N$ on log-log scale with $N = 2^{\bf s}, \,{\bf s} = 3,4,5,6$ at time $T=1$ for multiplicative noise with $g(u)=\sin(u).$ 
We simulate the exact solution with the numerical one using a sufficiently small step size $\Delta t = 2^{-11}.$  
It can be observed that the slopes of six fully-discrete schemes are closed to 1 on the temporal convergence order.

%---------------------------------------------------------------------------------------------------------------
\subsection{Nonlinear stochastic wave equation with $f(u) = u^3$} 
\label{sec;ne3} 
In this subsection, we focus on  the following nonlinear stochastic equation
\begin{equation} \label{eq6.2}\left\{\begin{aligned}
&du = vdt,  \\
&dv = \Delta u dt - u^3 dt +g(u) dW(t), &&(x, t) \in (0, 1) \times (0,100],
\end{aligned} \right. \end{equation}
subject to the initial conditions $u(x,0)=0,$ $v(x,0)=1$. The diffusion coefficients are the same as those in Section \ref{sec;ne1}, and the six fully-discrete schemes are chosen as same as those in Section \ref{sec;ne2}.  Note that the well-posedness of the nonlinear stochastic wave equation  \eqref{eq6.2}  with cubic nonlinearity and the considered diffusion coefficient can be obtained in \cite{Barbu, Chow, Schurz}.

\begin{figure}
\centering
\subfigure{
\begin{minipage}{6cm}
\centering
\includegraphics[height=5cm,width=6cm]{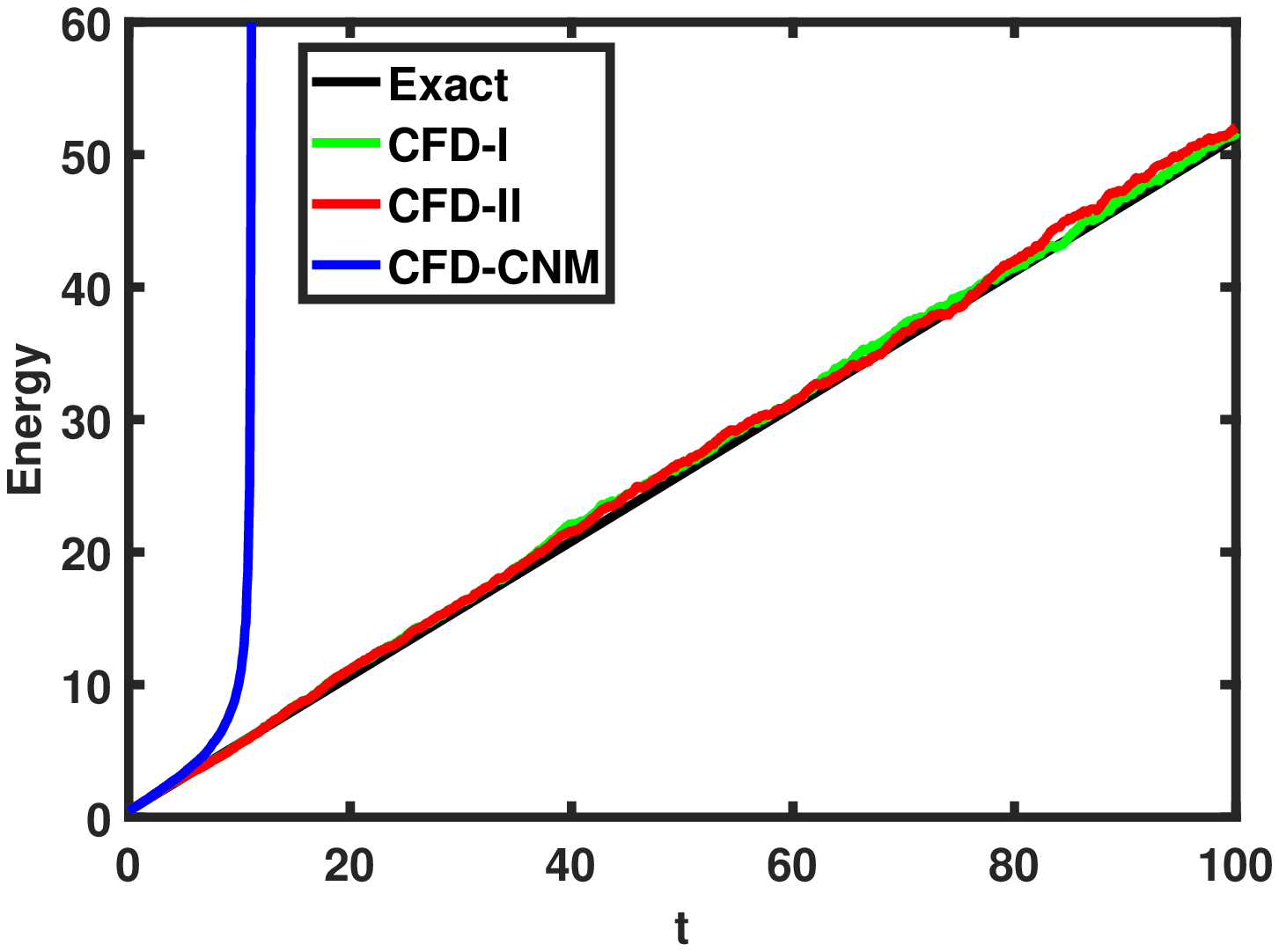}
\end{minipage}
}
\subfigure{
\begin{minipage}{6cm}
\centering
\includegraphics[height=5cm,width=6cm]{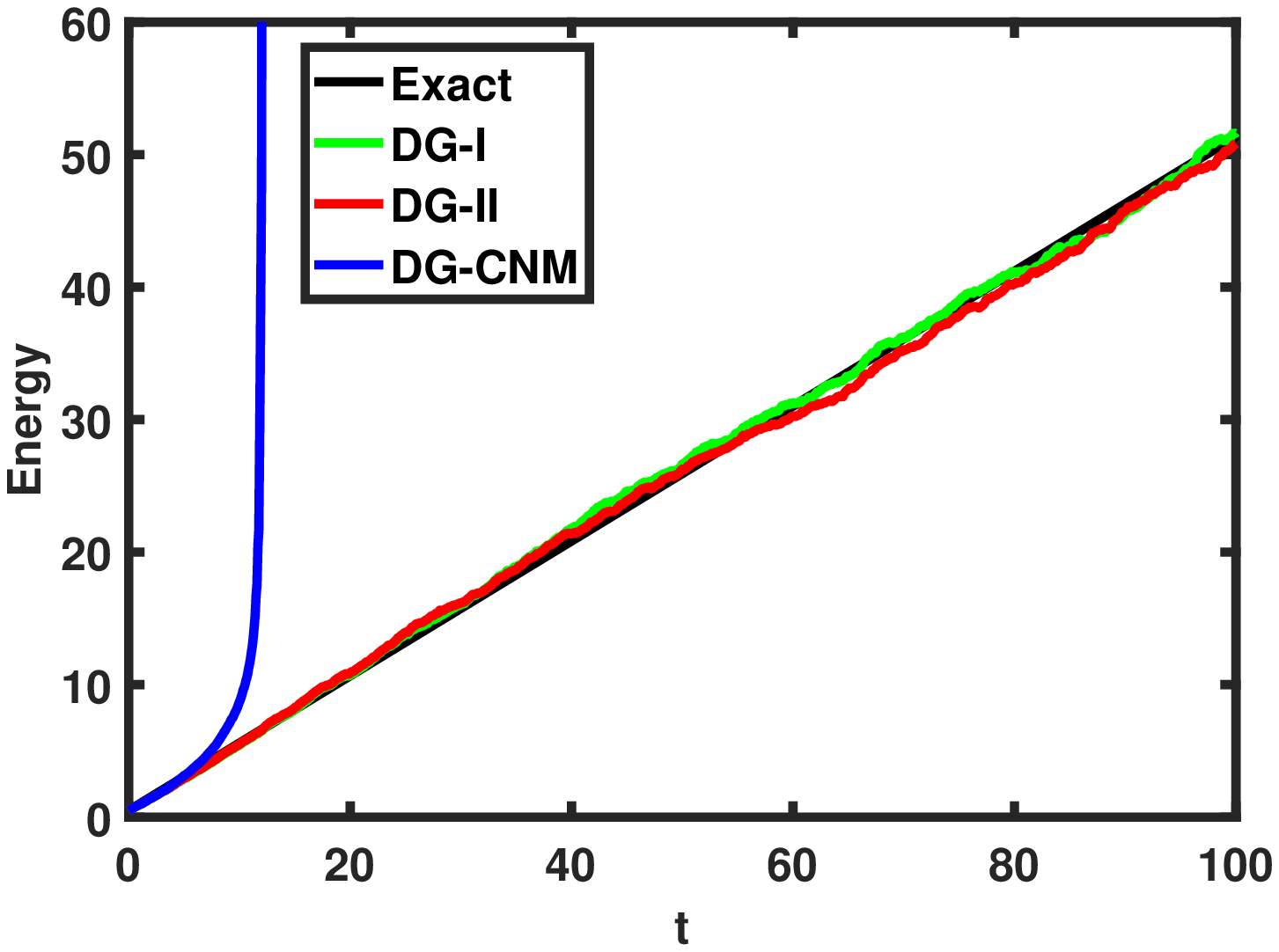}
\end{minipage}
}
\caption{ Averaged energy evolution relationship ($f(u) = u^3, g(u) = 1$) with $\Delta t= 1/12, h=1/20$}
\label{fig6}
\end{figure}

\begin{figure}
\centering
\subfigure{
\begin{minipage}{6cm}
\centering
\includegraphics[height=5cm,width=6cm]{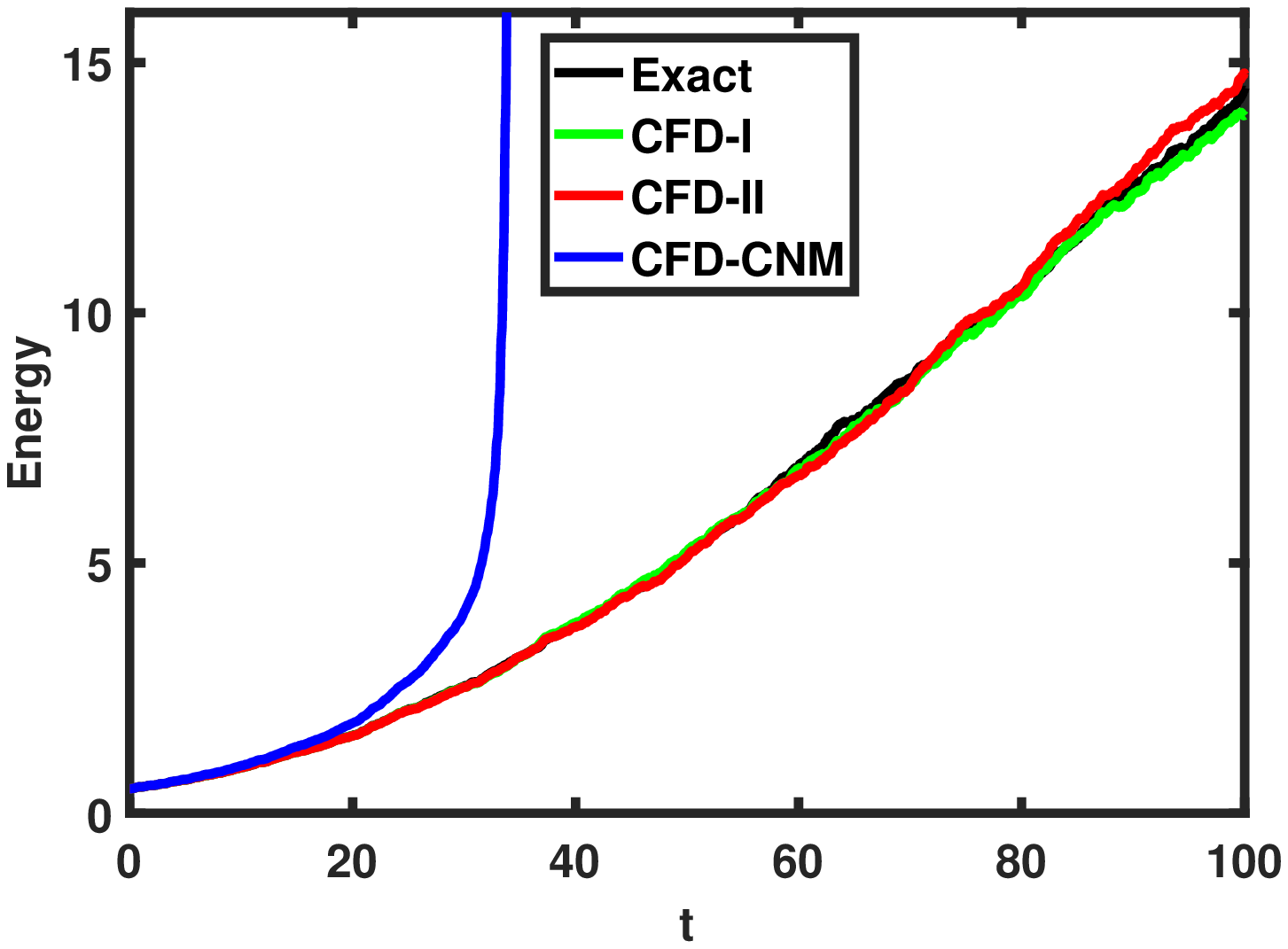}
\end{minipage}
}
\subfigure{
\begin{minipage}{6cm}
\centering
\includegraphics[height=5cm,width=6cm]{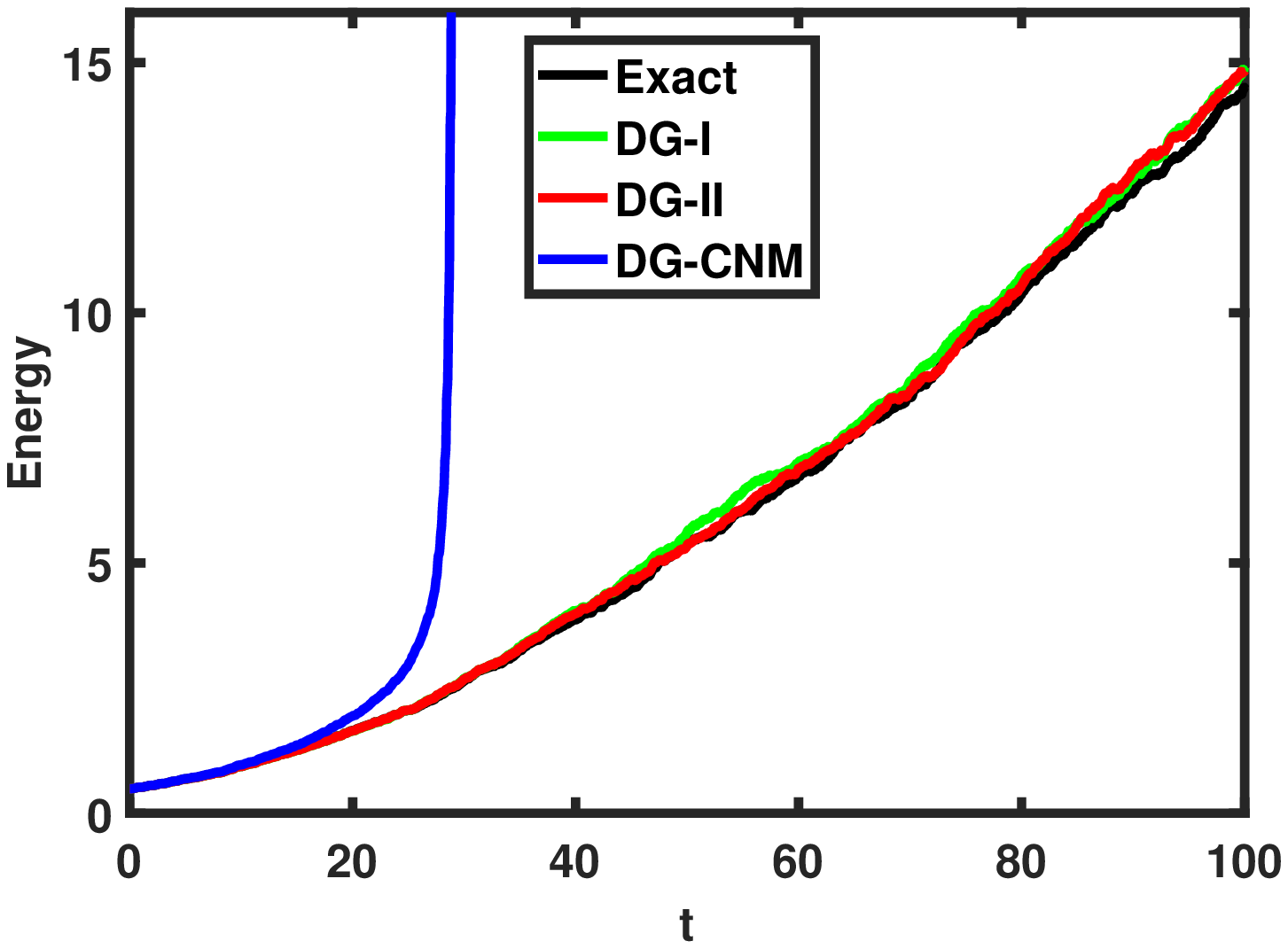}
\end{minipage}
}
\caption{ Averaged energy evolution relationship ($f(u) = u^3, g(u)=\sin (u)$) with $\Delta t= 1/25, h=1/50$}
\label{fig7}
\end{figure}

\begin{figure}
\centering
\subfigure{
\begin{minipage}{6cm}
\centering
\includegraphics[height=5cm,width=6cm]{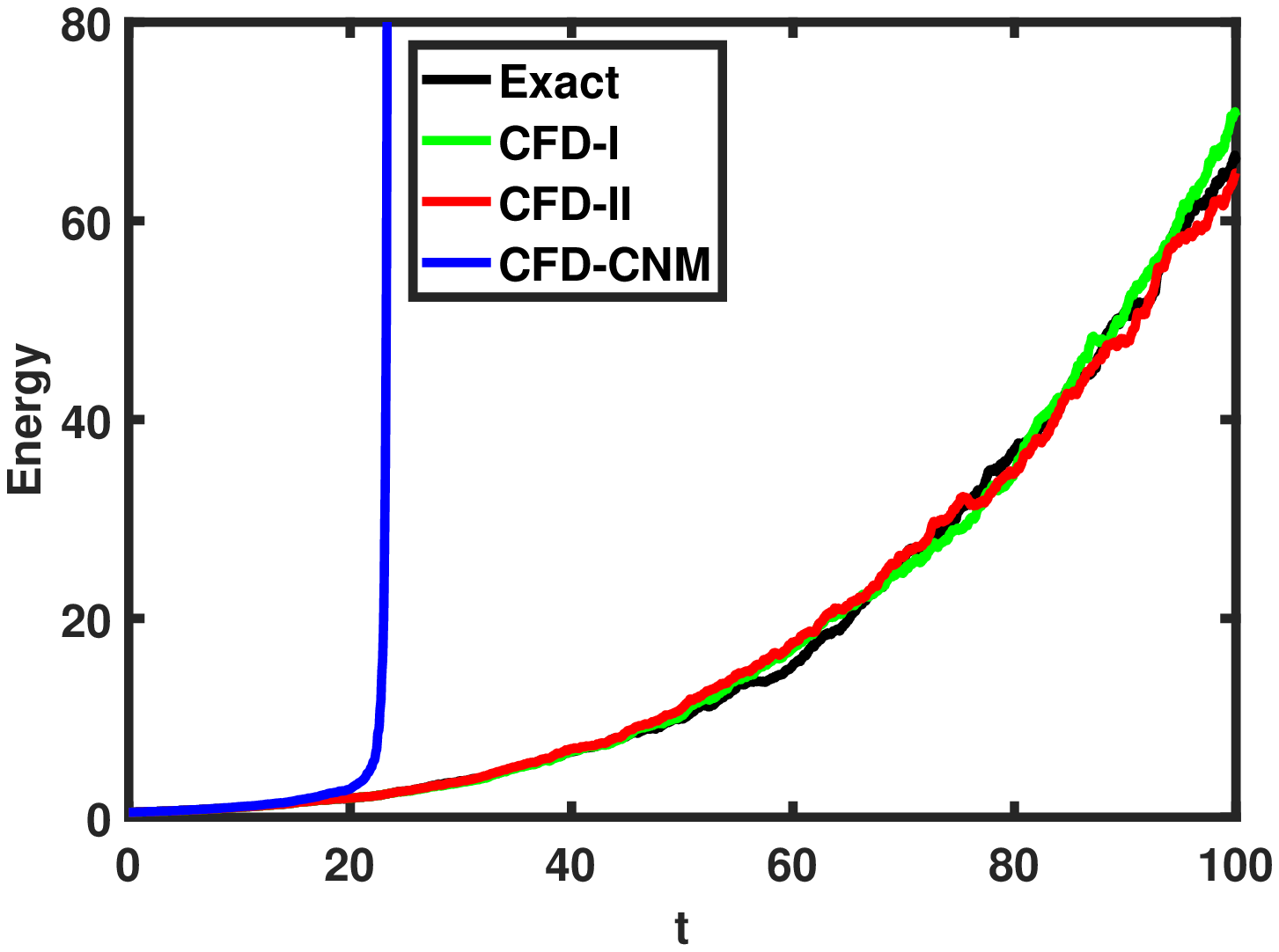}
\end{minipage}
}
\subfigure{
\begin{minipage}{6cm}
\centering
\includegraphics[height=5cm,width=6cm]{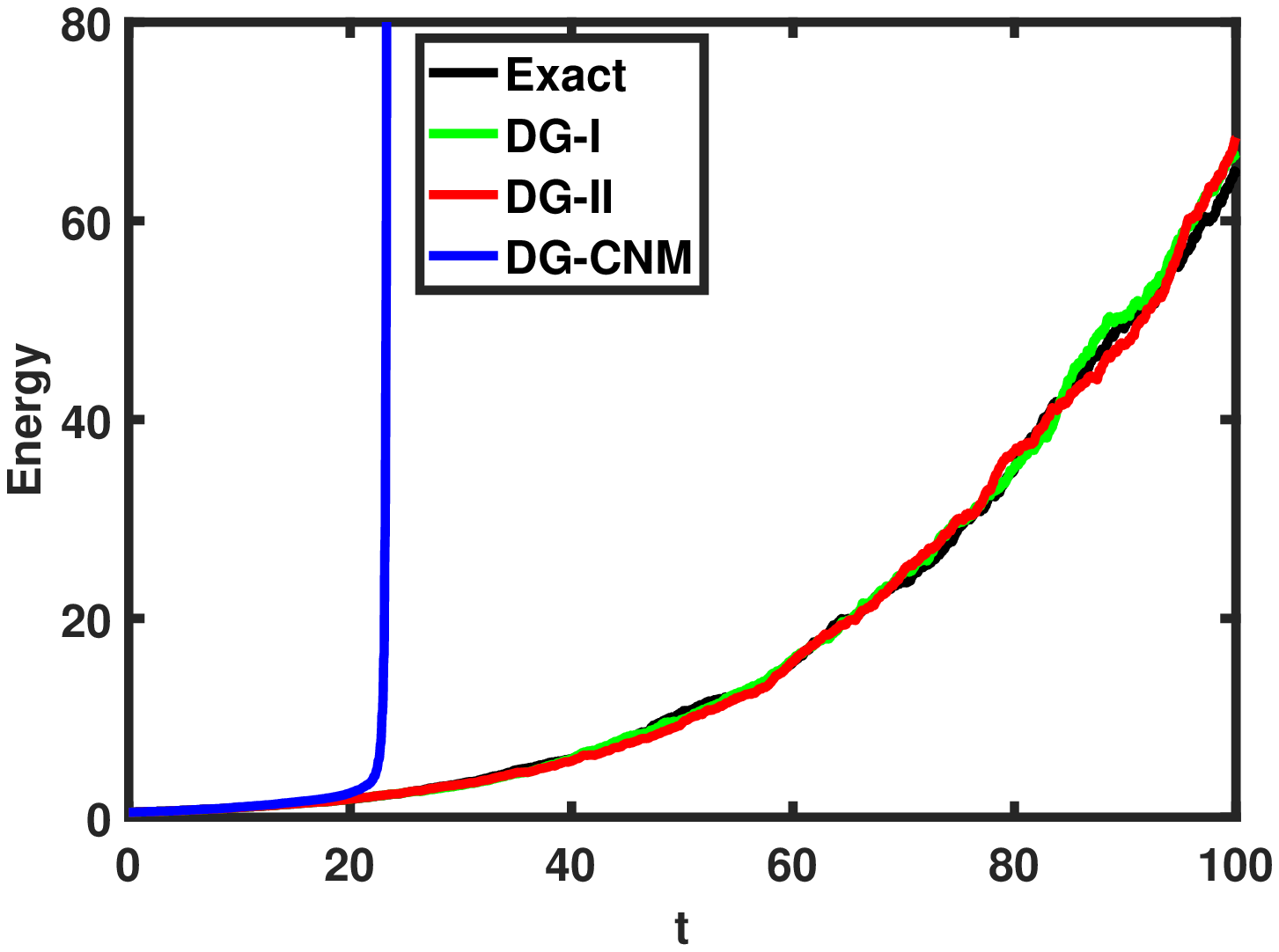}
\end{minipage}
}
\caption{ Averaged energy evolution relationship ($f(u) = u^3, g(u)=u$) with $\Delta t= 1/20, h=1/50$}
\label{fig8}
\end{figure}

For the case of additive noise,  Fig. \ref{fig6} shows that the averaged energies of CFD-I, CFD-II, DG-I and DG-II grow linearly, and coincide with the averaged energy of the exact solution.  However, the averaged energies of CFD-CNM and DG-CNM grow very fast and do not grow linearly, which fail to preserve the averaged energy evolution law of the original system.
 For the cases of multiplicative noise, i.e., the diffusion coefficients are chosen as $g(u) = \sin(u)$ and $g(u) = u$, Figs. \ref{fig7}-\ref{fig8} also indicate that the four proposed schemes preserve the evolution of averaged energy,  while the CFD-CNM and DG-CNM schemes do not, which coincide with theoretical analysis results,  and Fig. \ref{fig8} demonstrates clearly the effectiveness of the proposed fully-discrete schemes for simulating the nonlinear stochastic wave equation with multiplicative noise and non-globally Lipschitz continuous drift coefficient. In addition, it can be seen clearly that the averaged energy in Fig. \ref{fig8} grows faster than that in Fig. \ref{fig7}. Compared with the reference line in Fig. \ref{fig6.5} , it can be observed that the strong convergence order of six fully-discrete schemes is 1 in  temporal direction.

\iffalse
{\color{red}
	Furthermore, when the drift coefficient $f(u)$ and diffusion coefficient $g(u)$ in  \eqref{eq1.1} are nonlinear functions  satisfying that  $f(u)$ is a polynomial of the form $f(u) = \sum_{i=1}^{2k+1}u^i$, where $k\geq 0$, and there exist positive constants $\alpha_1 \geq 0$ and $\alpha_2 >0$ such that  $-\sum_{i=1}^{2k+1}\frac{2}{i+1}u^{i+1}\geq (\alpha_1 + \alpha_2 u^{2k})u^2$ for $\forall~u\in \mathscr H$. $g(u)$ is continuous, and there exist positive constants $C_1$ and $C_2$ such that 
	\begin{equation*}
	\| g(u)\|_{\mathscr{L}_2^0} ^2 \leq C_1(1+\|u\|^{2(k+1)}), ~ \|g(u_1)-g(u_2)\|^2_{\mathscr{L}_2^0}
	\leq C_2(1+\|u_1\|^{2k} +\|u_2\|^{2k}\|u_1-u_2\|^2 ),~\forall~ u, u_1, u_2\in\mathscr H, 
	\end{equation*}
	then \eqref{eq1.1} admits a  unique continuous global solution. The well-posedness of problem \eqref{eq1.1},  under the case of  high dimension and  more general drift  and diffusion coefficients , can refer . 
}
\fi

%----------总结--------------------------------------------------------------------------------------------------
\section{ Conclusion}
In this paper, we propose fully-discrete schemes by the compact finite difference method or the interior penalty discontinuous Galerkin finite element method in space, the discrete gradient method and the Pad\'e approximation in time, for solving the nonlinear stochastic wave equations driven by multiplicative noise.  We prove that the proposed schemes preserve the discrete averaged energy evolution laws exactly. Numerical experiments confirm the theoretical analysis results. One future work is the study of the strong convergence analysis and the estimate of the strong convergence order for the proposed schemes,  which is widely used to characterize the efficiency and accuracy of numerical method.  Another future work is the numerical study of the nonlinear stochastic wave equations with non-globally Lipschitz nonlinearity, which is difficult to obtain the well-posedness and strong convergence order of numerical schemes.

\section*{Acknowledgements}
This work is supported by National Natural Science Foundation of China (No. 11971470, No. 11871068,
No. 12031020, No. 12022118).

%----------参考文献-----------------------------------------------------------------------------

\end{document}